\documentclass{amsart}

\usepackage{epsfig}
\usepackage{amscd}
\usepackage{amsmath, amssymb, amsthm}
\usepackage{latexsym}
\usepackage{graphics}
\usepackage{color}

\newtheorem*{main-theorem}{Main Theorem}
\newtheorem{proposition}{Proposition}[section]
\newtheorem{theorem}{Theorem}
\newtheorem*{old-thm}{Theorem}
\newtheorem{lemma}[proposition]{Lemma}
\newtheorem{corollary}[theorem]{Corollary}

\theoremstyle{definition}

\newtheorem*{remark}{Remark}
\newtheorem{claim}[proposition]{Claim}

\numberwithin{equation}{section}

\def\reals{{\mathbb R}}
\def\cx{{\mathbb C}}

\def\Ci{{\mathcal C}^\infty}
\def\Re{\,\mathrm{Re}\,}
\def\Im{\,\mathrm{Im}\,}
\def\tom{\widetilde{\omega}}

\def\WF{\mathrm{WF}_h\,}

\def\supp{\mathrm{supp}\,}

\def\id{\,\mathrm{id}\,}

\def\O{{\mathcal O}}
\def\SS{{\mathbb S}}
\def\s{{\mathcal S}}

\def\nbhd{\mathrm{neigh}\,}
\def\neigh{\mathrm{neigh}\,}
\def\Op{\mathrm{Op}\,}

\def\esssupp{\text{ess-supp}\,}
\def\ad{\mathrm{ad}\,}

\def\hsc{{\left( \tilde{h}/h \right)}}

\def\csh{{\left( h/\tilde{h} \right)}}
\def\phi{\varphi}
\def\half{{\frac{1}{2}}}
\def\dist{\text{dist}\,}
\def\diag{\text{diag}\,}
\def\l{{L^2(X)}}
\def\be{\begin{eqnarray*}}
\def\ben{\begin{eqnarray}}
\def\ee{\end{eqnarray*}}
\def\een{\end{eqnarray}}
\def\lll{\left\langle}
\def\rrr{\right\rangle}
\def\contraction{\Bigg\lrcorner}

\begin{document}
\title[Non-concentration]{Semiclassical Non-concentration near
  Hyperbolic Orbits}
\author{Hans Christianson}
\address{Department of Mathematics, University of California, Berkeley, CA 94720 USA}
\email{hans@math.berkeley.edu}

\begin{abstract}
For a large class of semiclassical pseudodifferential operators, including Schr\"odinger operators,
$ P ( h ) = -h^2 \Delta_g + V (x) $, on compact Riemannian manifolds, we give logarithmic lower bounds on 
the mass of eigenfunctions outside neighbourhoods of generic closed hyperbolic orbits. More precisely we show that
if $ A $ is a pseudodifferential operator which is microlocally equal to the identity near the 
hyperbolic orbit and microlocally zero away from the orbit, then
\[  \| u \| \leq C (\sqrt{\log(1/h)}/ h ) \| P (h)u \| + C \sqrt {\log(1/h )} \| ( I - A ) u \| \,. \]
This generalizes earlier estimates of Colin de Verdi\`ere-Parisse \cite{CVP} obtained for a special case, and
of Burq-Zworski \cite{BZ} for real hyperbolic orbits. 
\end{abstract}
\maketitle


\section{Introduction}
\label{intro}
\indent 
To motivate the general result, we first present two applications.  If
$(X,g)$ is a Riemannian manifold with Laplacian $\Delta_g$, we consider the eigenvalue problem
\be
-\Delta_g u = \lambda^2 u, \ \ \ 
\|u\|_{L^2(X)} = 1. 
\ee
If $U$ is a small neighbourhood of a closed {\em hyperbolic} geodesic $\gamma$, 
we show that 
\be
\int_{X \setminus U} |u|^2 dx \geq \frac{c}{ \log | \lambda
    |},
\ee
that is, if $u$ 
concentrates near $\gamma$, the rate is logarithmic. This generalizes
results of Colin de Verdi\`ere-Parisse \cite{CVP} and Burq-Zworski
\cite{BZ}.

As another application of our main results we consider the damped wave equation
\begin{eqnarray*}
\left\{ \begin{array}{l}
\left( \partial_t^2 - \Delta + 2a(x) \partial_t \right) u(x,t) = 0, \quad (x,t) \in X \times (0, \infty) \\
u(x,0) = 0, \quad \partial_t u(x,0) = f(x).
\end{array} \right.
\end{eqnarray*}
We prove in \S \ref{wave-section} that if $a(x)>0$ outside a
neighbourhood of a closed hyperbolic geodesic $\gamma$, we have the following energy estimate:
\be
\| \partial_t u \|_{L^2(X)}^2 + \left\| \nabla u
\right\|_{L^2(X)}^2 \leq C e^{-t/C} \|f\|_{H^\epsilon(X)}^2,
\ee
for all $\epsilon>0$.  (In \S \ref{wave-section} a weaker
geometric control condition in the spirit of Rauch-Taylor \cite{RT} is
considered.) This application was suggested to us by M. Hitrik, 
and it generalizes an example of Lebeau \cite{Leb}.

We now turn to the general case. Let $X$ be a compact $n$-dimensional manifold without boundary.  We consider a selfadjoint 
pseudodifferential operator, $P(h)$, with real principal symbol $p$. 
We assume throughout if $p=0$ then $dp \neq 0$, and that $p$ is elliptic outside
of a compact subset of $T^*X$.
Assume that 
$$ \gamma \subset  p^{-1} ( 0 ) $$ 
is  a closed loxodromic orbit of the Hamiltonian flow of $p$.  
Let $N \subset \{p=0\}$ be a Poincar\'{e} section for $\gamma$ and let $S$ be the Poincar\'e map.  
The assumption that $\gamma$ be loxodromic means that no eigenvalue of $dS(0,0)$ lies on the unit circle.  
We assume also that $dS(0,0)$ has no real negative eigenvalues.

\begin{main-theorem}
Let $A \in \Psi_h^{0,0}$ be a pseudodifferential operator whose principal symbol is 
$1$ near $\gamma$ and $0$ away from $\gamma$. Then, there exist constants $ h_0>0
$ and $0 < C < \infty$ so that we have uniformly in $ 0 < h < h_0 $,
\begin{eqnarray}
\label{main-theorem-3-est}
\| u \| \leq C \frac{ \sqrt{\log ( 1/h )}} h \| P ( h ) u \| + C
\sqrt{ \log (1/h ) } \| ( I - A ) u \| \,, 
\end{eqnarray}
where the norms are $ L^2 $ norms on $ X$.
In particular if a family,  $u = u ( h ) $ satisfies
\begin{eqnarray*}
P(h) u = \O_{L^2}(h^\infty),  \ \ \ 
\| u \|_{L^2(X)} = 1\, ,
\end{eqnarray*} 
then
\begin{eqnarray}
\label{main-theorem-3-est-2}
 \left\| (I - A) u \right\|_{L^2(X)} \geq \frac{1}{C} \log \left(\left(1/h\right)\right)^{-\half} \,, \ \ 0 < h < h_0 \,.
\end{eqnarray}
\end{main-theorem}
We note that the assumptions on $ A $ imply that 
$\WF (A) $ is contained in a neighbourhood of $\gamma$, while $\WF(I - A)$ is away from $\gamma$, see 
\S \ref{preliminaries} for definitions.

Colin de Verdi\`ere and Parisse \cite{CVP} have shown that the
estimates (\ref{main-theorem-3-est}-\ref{main-theorem-3-est-2}) are sharp in the case where 
$ X $ is a segment of a hyperbolic cylinder and $ P ( h ) = -h^2 \Delta_g $ is its Dirichlet Laplacian.
Even though the closed orbit at the ``neck'' of the cylinder is hyperbolic, the flow is completely 
integrable in that case. This shows that eliminating the  $\log (h^{-1})$ factor requires global 
conditions on the classical flow.

The assumption that the Poincar{\'e} map has no negative eigenvalues
is standard in the literature on quantum Birkhoff normal forms (see, for example, \cite{IaSj},
\cite{ISZ}, and \cite{Ze}), and in the present work serves to
eliminate cases in which current techniques seem to break down.  It is
important to note that this case does arise,
as in the example in \cite{Kl} \S $3.4$.

There are many examples in which the hypotheses of the theorem are
satisfied, the simplest of which is the case in which $p = |\xi|^2- E(h)$ for $E(h)>0$.
Then the Hamiltonian flow of $p$ is the geodesic flow, so if the
geodesic flow has a closed hyperbolic orbit, there is
non-concentration of eigenfunctions, $u(h)$, for the equation
\be
-h^2 \Delta u(h) = E(h) u(h).
\ee
Another example of such a $p$ is the case $p = |\xi|^2 + V(x)$, where
$V(x)$ is a confining potential with three ``bumps'' or ``obstacles''
in the lowest energy level (see Figure \ref{fig:fig13}).  In the
appendix to \cite{Sjo2a} it is shown that for an interval of energies
$V(x) \sim 0$, there is a closed hyperblic orbit $\gamma$ of the Hamiltonian
flow which ``reflects'' off the bumps (see Figure \ref{fig:fig14}).
Loxodromic orbits may be constructed by considering $3$-dimensional
hyperbolic billiard problems (see, for example, \cite{AuMa}), although
in the present work we are assuming the orbit does not intersect
the boundary of the manifold.  In addition, Proposition
\ref{normal-prop-1} gives a somewhat artificial means of constructing a manifold diffeomorphic to a neighbourhood in $T^*\SS_{(t, \tau)}^1
\times T^* \reals_{(x, \xi)}^{n-1}$ which contains a loxodromic orbit $\gamma$
by starting with the Poincar\'{e} map $\gamma$ is to have.


\begin{figure}
\centerline{\epsfig{figure=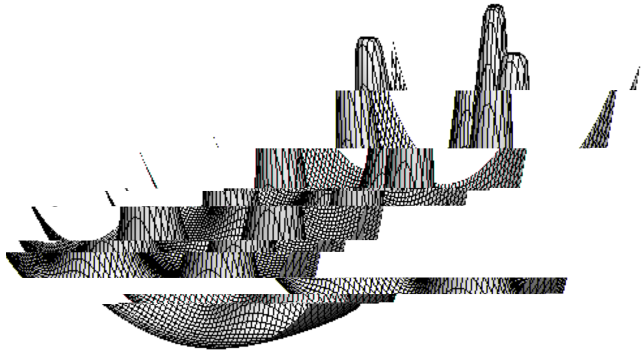, height=2.5in}}
\caption{\label{fig:fig13} A confining potential $V(x)$ with three
  bumps at the lowest energy level $E <0$.}
\end{figure}

\begin{figure}
\include{fig14}
\caption{\label{fig:fig14} The level set $V(x) = 0$ and the closed
  hyperbolic orbit $\gamma$.}
\end{figure}

\indent In order to prove the Main Theorem, we will first prove that 
the principal symbol of $P(h)$ can be put into a normal form near $\gamma$. This 
will allow analysis of  small complex perturbations of $P(h)$.  These are defined as follows: 
let $a \in \Ci(T^*X,[0,1])$ be equal to $0$ in a neighbourhood of $\gamma$ and $1$ outside of a larger 
neighbourhood of $ \gamma $.  For $z \in [-1,1] + i[-\delta, \delta]$, define
\begin{eqnarray}
\label{Q(z)}
Q(z):= P(h) - z - ih Ca^w,
\end{eqnarray}
for a constant $C$ to be fixed later.  The following theorem states that by perturbing $P(h)$ into $Q(z)$
 we are able to push the spectrum of $P(h)$ into the lower half plane.
\begin{theorem}
\label{main-theorem-1}
There exist constants $c_0>0$, $h_0>0$, and $N_0$ such that for $u$ with $\WF(u)$ in a sufficiently small
 neighbourhood of $\gamma$, $z \in [-1,1] + i(-c_0h, + \infty)$, and $0 < h < h_0$ we have
\begin{eqnarray}
\label{main-theorem-1-est}
\left\| Q(z) u \right\|_\l \geq C^{-1} h^{N_0} \left\| u \right\|_\l
\end{eqnarray}
for some constant $C$.
\end{theorem}
Using Theorem \ref{main-theorem-1} and a semiclassical adaptation of the ``three-lines'' theorem from complex analysis, 
we will be able to deduce the following estimate.
\begin{theorem}
\label{main-theorem-2}
Suppose $Q(z)$ is given by \eqref{Q(z)}, and $z \in I \Subset (-\infty, \infty)$.  Then there 
is $h_0 >0$ and $0 < C < \infty$ such that for $0 < h < h_0 $,
\begin{eqnarray}
\label{main-theorem-2-est-1}
\left\| Q(z)^{-1} \right\|_{L^2(X) \to L^2(X)} \leq C \frac{\log (1/h) }{h}.
\end{eqnarray}
If $\phi \in \Ci_c(X)$ is supported away from $\gamma$, then
\begin{eqnarray}
\label{main-theorem-2-est-2}
\left\| Q(z)^{-1} \phi \right\|_{L^2(X) \to L^2(X)} \leq C \frac{\sqrt{ \log (1/h)}}{h}.
\end{eqnarray}
\end{theorem}
In order to apply the results of Theorems \ref{main-theorem-1} and \ref{main-theorem-2} to the Main Theorem, 
we observe that for $A$ as in the statement of the Main Theorem we
have $Q(0)A = P(h)A$ microlocally and apply a commutator argument.

This note is organized as follows.  \S \ref{preliminaries} recalls basic facts about 
 $h$-pseudodifferential operators on manifolds.  This is followed in \S \ref{FIO} with a review 
of some standard results from the theory of $h$-Fourier Integral Operators.  In \S \ref{symplectic-geometry} we present 
some symplectic geometry and prove the principal symbol can be put
 into a normal form in the case all the eigenvalues of $dS(0)$ are distinct.  
\S \ref{main-theorem-1-proof} contains the proof of Theorem
 \ref{main-theorem-1} in the case of distinct eigenvalues, then re-examines the normal
 form of the principal symbol to show how it may be extended to the
 case when the eigenvalues are not distinct, and contains the details
 of the more general case of Theorem \ref{main-theorem-1}.  Finally, in \S \ref{main-theorem-2-proof} we prove Theorem
 \ref{main-theorem-2} and the Main Theorem.  In \S
 \ref{wave-section} we follow a suggestion of M. Hitrik to apply the
 techniques of \S \ref{symplectic-geometry}-\ref{main-theorem-2-proof}
 to the damped wave equation. 

The impetus for this paper came when M. Zworski suggested generalizing results from the appendix of \cite{BZ}, 
as well as correcting a mistake which was discovered by J.-F. Bony, S. Fujiie, T. Ramond, 
and M. Zerzeri (see \cite{BFRZ} for their closely related work).  This paper generalizes the statements of the theorems 
from the case of real hyperbolic trajectories to complex hyperbolic or loxodromic trajectories as well as 
correcting the mistake. 

{\bf Acknowledgements:} The author would like to thank Maciej Zworski for much help
  and support during the writing of this paper, Alan Weinstein and Fr\'ed\'eric Naud for helpful
  conversations, as well as the NSF for partial support.  He would
  like to thank Michael Hitrik for comments on an early draft and suggesting Section
  \ref{wave-section}, and Laurent Thomann and Steve Zelditch for careful reading of an
  early draft of this paper.


\section{Preliminaries}
\label{preliminaries}
This section contains some basic definitions and results from semiclassical and microlocal 
analysis which we will be using throughout the paper.  This is essentially standard, but we 
include it for completeness.  We will follow the presentation in \cite{BZ}, \S $2$.  
Let $X$ be a smooth, compact manifold.  We will be operating on half-densities,
\begin{eqnarray*}
u(x)|dx|^{\frac{1}{2}} \in \Ci\left(X, \Omega_X^{\frac{1}{2}}\right),
\end{eqnarray*}
with the informal change of variables formula
\begin{eqnarray*}
u(x)|dx|^{\frac{1}{2}} = v(y)|dy|^\half, \,\, \text{for}\,\, y = \kappa(x) 
\Leftrightarrow v(\kappa(x))|\kappa'(x)|^\half = u(x).
\end{eqnarray*}
By symbols on $X$ we mean
\begin{eqnarray*}
\lefteqn{\s^{k,m} \left(T^*X, \Omega_{T^*X}^\half \right):= } \\
& = &\left\{ a \in \Ci(T^*X \times (0,1], \Omega_{T^*X}^\half): \left| \partial_x^\alpha 
\partial_\xi^\beta a(x, \xi; h) \right| \leq C_{\alpha \beta}h^{-m} \langle \xi \rangle^{k - |\beta|} \right\}.
\end{eqnarray*}
There is a corresponding class of pseudodifferential operators $\Psi_h^{k,m}(X, \Omega_X^\half)$ 
acting on half-densities defined by the local formula (Weyl calculus) in $\reals^n$:
\begin{eqnarray*}
\Op_h^w(a)u(x) = \frac{1}{(2 \pi h)^n} \int \int a \left( \frac{x + y}{2}, \xi; h \right) 
e^{i \langle x-y, \xi \rangle / h }u(y) dy d\xi.
\end{eqnarray*}
We will occasionally use the shorthand notations $a^w := \Op_h^w(a)$ and $A:=\Op_h^w(a)$ when 
there is no ambiguity in doing so. \\
\indent We have the principal symbol map
\begin{eqnarray*}
\sigma_h : \Psi_h^{k,m} \left( X, \Omega_X^\half \right) \to \s^{k,m} \left/ \s^{k, m-1} 
\left(T^*X, \Omega_{T^*X}^\half \right) \right.,
\end{eqnarray*}
which gives the left inverse of $\Op_h^w$ in the sense that 
\begin{eqnarray*}
\sigma_h \circ \Op_h^w: \s^{k,m} \to \s^{k,m}/\s^{k, m-1} 
\end{eqnarray*}
is the natural projection.  Acting on half-densities in the Weyl calculus, the principal 
symbol is actually well-defined in $\s^{k,m} / \s^{k, m-2}$, that is, up to $\O(h^2)$ in 
$h$ (see, for example \cite{EvZw} Appendix D).  

We will use the notion of wave front sets for pseudodifferential operators on manifolds.  
If $a \in \s^{k,m}(T^*X, \Omega_{T^*X}^\half)$, we define the singular support or essential support for $a$:
\begin{eqnarray*}
\esssupp_h a \subset T^*X \bigsqcup \SS^*X,
\end{eqnarray*}
where $\SS^*X = (T^*X \setminus \{0\}) / \reals_+$ is the cosphere bundle (quotient taken with 
respect to the usual multiplication in the fibers), and the union is disjoint.  $\esssupp_h a$ is defined using complements:
\begin{eqnarray*}
\lefteqn{\esssupp_h a := } \\
& = & \complement \left\{ (x, \xi) \in T^*X : \exists \epsilon >0, \,\,\, \partial_x^\alpha 
\partial_\xi^\beta a(x', \xi') = \O(h^\infty), \,\,\, d(x, x') + |\xi - \xi'| < \epsilon \right\} \\
&& \bigcup \complement \{ (x, \xi) \in T^*X \setminus 0 : \exists \epsilon > 0, \,\,\, \partial_x^\alpha 
\partial_\xi^\beta a(x', \xi') = \O (h^\infty \langle \xi \rangle^{-\infty}),  \\
&& \quad \quad \quad d(x, x') + 1 / |\xi'| + | \xi/ |\xi| - \xi' /
|\xi'| | < \epsilon \} / \reals_+.
\end{eqnarray*}
We then define the wave front set of a pseudodifferential operator $A \in \Psi_h^{k,m}( X, \Omega_X^\half )$:
\begin{eqnarray*}
\WF(A) : = \esssupp_h(a), \,\,\, \text{for} \,\,\, A = \Op_h^w(a).
\end{eqnarray*}
Finally for distributional half-densities $u \in \Ci( (0, 1]_h, \mathcal{D}'(X, \Omega_X^\half))$ 
such that there is $N_0$ so that $h^{N_0}u$ is bounded in $\mathcal{D}'(X, \Omega_X^\half)$, we can 
define the semiclassical wave front set of $u$, again by complement:
\begin{eqnarray*}
\lefteqn{\WF (u) := } \\
&= & \complement \{(x, \xi) : \exists A \in \Psi_h^{0,0}, \,\, \text{with} \,\, \sigma_h(A)(x,\xi) \neq 0, \,\,  \\ 
&& \quad \text{and} \,\,Au \in h^\infty \Ci((0,1]_h, \Ci(X, \Omega_X^\half)) \}.
\end{eqnarray*}
For $A = \Op_h^w(a)$ and $B = \Op_h^w(b)$, $a \in \s^{k,m}$, $b \in \s^{k',m'}$ we have the composition 
formula (see, for example, \cite{DiSj})
\begin{eqnarray}
\label{Weyl-comp}
A \circ B = \Op_h^w \left( a \# b \right),
\end{eqnarray}
where
\begin{eqnarray}
\label{a-pound-b}
\s^{k + k', m+m'} \ni a \# b (x, \xi) := \left. e^{\frac{ih}{2} \omega(Dx, D_\xi; D_y, D_\eta)} 
\left( a(x, \xi) b(y, \eta) \right) \right|_{{x = y} \atop {\xi = \eta}} ,
\end{eqnarray}
with $\omega$ the standard symplectic form. \\
\indent We will need the definition of microlocal equivalence of operators.  Suppose 
$T: \Ci(X, \Omega_X^\half) \to \Ci(X, \Omega_X^\half)$ and that for any seminorm $\| \cdot \|_1$ on 
$\Ci(X, \Omega_X^\half)$ there is a second seminorm $\| \cdot \|_2$ on $\Ci(X, \Omega_X^\half)$ such that 
\begin{eqnarray*}
\| Tu\|_1 = \O(h^{-M_0})\|u \|_2
\end{eqnarray*}
for some $M_0$ fixed.  Then we say $T$ is {\it semiclassically tempered}.  We assume for the rest of 
this paper that all operators satisfy this condition.  Let $U,V \subset T^*X$ be open precompact sets.  
We think of operators defined microlocally near $V \times U$ as equivalence classes of tempered operators.  
The equivalence relation is
\begin{eqnarray*}
T \sim T' \Longleftrightarrow A(T-T')B = \O(h^\infty): \mathcal{D}'\left( X, \Omega_X^\half \right) \to \Ci 
\left(X, \Omega_X^\half \right)
\end{eqnarray*}
for any $A,B \in \Psi_h^{0,0}(X, \Omega_X^\half)$ such that 
\begin{eqnarray*}
&& \WF (A) \subset \widetilde{V}, \quad \WF (B) \subset \widetilde{U}, \,\, \text{with} \,\, \widetilde{V}, 
\widetilde{U} \,\, \text{open and } \\
&& \quad \quad \overline{V} \Subset \widetilde{V} \Subset T^*X, \quad \overline{U} \Subset 
\widetilde{U} \Subset T^*X.
\end{eqnarray*}
In the course of this paper, when we say $P=Q$ {\it microlocally} near $U \times V$, we mean for any $A$, $B$ as above,
\begin{eqnarray*}
APB - AQB = \O_{L^2 \to L^2}\left( h^\infty \right),
\end{eqnarray*}
or in any other norm by the assumed precompactness of $U$ and $V$.  Similarly, we say $B = T^{-1}$ on $V \times V$ if 
$BT = I$ microlocally near $U \times U$ and $TB = I$ microlocally near $V \times U$. \\
\indent For this paper, we will need the following semiclassical version of Beals's Theorem (see \cite{DiSj} for a proof).  
Recall for operators $A$ and $B$, the notation $\ad_BA$ is defined as 
\begin{eqnarray*}
\ad_B A = \left[B, A \right].
\end{eqnarray*}
\begin{old-thm}[Beals's Theorem]
Let $A: \s \to \s'$ be a continuous linear operator.  Then $A = \Op_h^w(a)$ for a symbol $a \in \s^{0,0}$ if and only 
if for all $N \in \mathbb{N}$ and all linear symbols $l_1, \ldots l_N$, 
\begin{eqnarray*}
\ad_{\Op_h^w(l_1)} \circ \ad_{\Op_h^w(l_2)} \circ \cdots \circ \ad_{\Op_h^w(l_N)} A = \O(h^N)_{L^2 \to L^2}.
\end{eqnarray*}
\end{old-thm}
The following lemma (given more generally in \cite{BoCh}) will be used in the proof of Theorem \ref{main-theorem-1}.  
We include a sketch of the proof from \cite{SjZw2} here for completeness.  It is easiest to phrase in terms of 
order functions.  
A smooth function $m \in \Ci(T^*X; \reals)$ is called an order function if it satisfies
\begin{eqnarray*}
m(x, \xi) \leq C m(y, \eta) \left\langle \dist (x-y) + |\xi - \eta| \right\rangle^N
\end{eqnarray*}
for some $N \in \mathbb{N}$.  We say $a \in \s^l(m)$ if 
\begin{eqnarray*}
\left| \partial^\alpha  a \right| \leq C_{\alpha }h^{-l} m.
\end{eqnarray*}
If $l = 0$, we write $\s(m):= \s^0(m)$.
\begin{lemma}
\label{etG-lemma}
Let $m$ be an order function, and suppose $G \in \Ci (T^*X ; \reals)$ satisfies
\begin{eqnarray}
\label{G-cond-1}
G(x, \xi) - \log \left( m(x, \xi) \right) = \O (1),
\end{eqnarray}
and 
\begin{eqnarray}
\label{G-cond-2}
\partial_x^\alpha \partial_\xi^\beta G(x,\xi) = \O(1) \,\,\, \text{for} \,\,\, (\alpha, \beta) \neq (0,0).  
\end{eqnarray}
Then for $G^w = \Op_h^w(G)$ and $|t|$ sufficiently small,
\begin{eqnarray*}
\exp (tG^w) = \Op_h^w(b_t)
\end{eqnarray*}
for $b_t \in \s(m^t)$.  Here $e^{tG^w}$ is defined as the unique solution to the ordinary
differential equation
\be
\left\{ \begin{array}{l} \partial_t \left( U(t) \right) - G^w
  U(t) = 0 \\
U(0) = \id. \end{array} \right.
\ee
\end{lemma}
\begin{proof}[Sketch of Proof]
The conditions on $G$ \eqref{G-cond-1} and \eqref{G-cond-2} are equivalent to saying $e^{tG} \in \s(m^t)$.  We will 
compare $\exp tG^w$ and $\Op_h^w (\exp t G)$.
\begin{claim}
\label{U-inv-claim}
Set $U(t):= \Op_h^w (e^{tG}): \s \to \s$.  For $|t| < \epsilon_0$, $U(t)$ is invertible and $U(t)^{-1} = \Op_h^w (b_t)$ for 
$b_t \in \s(m^{-t})$, where $\epsilon_0$ depends only on $G$.
\end{claim}
\begin{proof}[Proof of Claim]
Using the composition law, we see $U(-t) U(t) = \id + \Op_h^w(E_t)$, with $E_t = \O(t)$.  Hence $\id + \Op_h^w(E_t)$ is 
invertible and using Beals's Theorem, we get $(\id + \Op_h^w(E_t))^{-1} = \Op_h^w(c_t)$ for $c_t \in \s(1)$.  
Thus $\Op_h^w(c_t)U(-t)U(t) = \id$, so
\begin{eqnarray*}
U(t)^{-1} = \Op_h^w \left( c_t \# \exp (-tG) \right),
\end{eqnarray*}
and subsequently $b_t \in \s(m^{-t})$.
\end{proof}
Now observe that
\begin{eqnarray*}
\frac{d}{dt} U(-t) = - \Op_h^w \left( G \exp(-tG) \right),  \,\,\, \text{and} \,\,\, U(-t) G^w = \Op_h^w \left( e^{-tG} \# G \right),
\end{eqnarray*}
so that
\begin{eqnarray}
\lefteqn{\frac{d}{dt} \left( U(-t) e^{tG^w} \right) =} \label{Op(A_t)}
\\ & = & -\Op_h^w \left( G \exp(-tG) \right) e^{tG^w} + 
\Op_h^w \left( e^{-tG} \# G \right) e^{tG^w} \nonumber \\
& = & \Op_h^w (A_t) e^{tG^w}, \nonumber
\end{eqnarray}
for $A_t \in \s(m^{-t})$.  To see \eqref{Op(A_t)}, recall that by the composition law,
\begin{eqnarray*}
e^{-tG} \# G = e^{-tG} G + \left( \text{terms with }\,G \,\,\, \text{derivatives} \right).
\end{eqnarray*}
Then the first terms in \eqref{Op(A_t)} will cancel and the remaining terms will all involve at least one 
derivative of $G$, which is then bounded by \eqref{G-cond-2}.  \\
\indent Set $C(t):= -\Op_h^w(A_t)U(-t)^{-1}$.  Claim \ref{U-inv-claim} implies $C(t) = \Op_h^w (c_t)$ for a family 
$c_t \in \s(1)$.  The composition law implies $c_t$ depends smoothly on $t$.  Then
\begin{eqnarray*}
\left( \frac{\partial}{\partial t} + C(t) \right) \left( U(-t) e^{tG^w} \right) = \Op_h^w (A_t) e^{tG^w} - 
\Op_h^w (A_t) e^{tG^w} = 0,
\end{eqnarray*}
so we have reduced the problem to proving the following claim.
\begin{claim}
Suppose $C(t) = \Op_h^w (c_t)$ with $c_t \in \s(1)$ depending smoothly on $t \in (-\epsilon_0, \epsilon_0)$.  If $Q(t)$ solves
\begin{eqnarray*}
\left\{ \begin{array}{c}
\left( \frac{\displaystyle \partial}{\displaystyle \partial t} + C(t) \right) Q(t) = 0, \\
Q(0) = \Op_h^w(q), \,\,\, \text{with} \,\,\, q \in \s(1),
\end{array} \right.
\end{eqnarray*}
then $Q(t) = \Op_h^w (q_t)$ with $q_t \in \s(1)$ depending smoothly on $t \in (-\epsilon_0, \epsilon_0)$.
\end{claim}
\begin{proof}[Proof of Claim]
The Picard existence theorem for ODEs implies $Q(t)$ exists and is bounded on $L^2$.  We want to use Beals's 
Theorem to show $Q(t)$ is actually a quantized family of symbols.  Let $l_1, \ldots, l_N$ be linear symbols.  
We will use induction to show that for any $N$ and any choice of the $l_j$, $\ad_{\Op_h^w (l_1)} \circ \cdots \circ 
\ad_{\Op_h^w(l_N)}Q(t) = \O(h^N)_{L^2\to L^2}$.  Since we are dealing with linear symbols, we take $h = 1$ for 
convenience.  First note
\begin{eqnarray*}
&&\frac{d}{dt} \ad_{\Op_h^w (l_1)} \circ \cdots \circ \ad_{\Op_h^w(l_N)} Q(t) +  \ad_{\Op_h^w (l_1)} \circ \cdots \circ 
\ad_{\Op_h^w(l_N)}\\
&& \quad \quad \quad \quad \quad \quad \quad \quad \quad \cdot  \left( C(t) Q(t) \right) = 0
\end{eqnarray*}
For the induction step, assume $ \ad_{\Op_h^w (l_1)} \circ \cdots \circ \ad_{\Op_h^w(l_k)}Q(t) = \O(1)$ is known 
for $k<N$ and observe
\begin{eqnarray*}
\lefteqn{ \ad_{\Op_h^w (l_1)} \circ \cdots \circ \ad_{\Op_h^w(l_N)} \left( C(t) Q(t) \right) = } \\
& = & C(t)  \ad_{\Op_h^w (l_1)} \circ \cdots \circ \ad_{\Op_h^w(l_N)} Q(t) + R(t),
\end{eqnarray*}
where $R(t)$ is a sum of terms of the form $A_k(t) \ad_{\Op_h^w (l_1)} \circ \cdots \circ \ad_{\Op_h^w(l_k)}Q(t)$ for 
each $k<N$ and $A_k(t) = \Op_h^w (a_k(t))$ with $a_k(t) \in \s(1)$.  Set $\tilde{Q}(t) =  \ad_{\Op_h^w (l_1)} \circ 
\cdots \circ \ad_{\Op_h^w(l_N)}Q(t)$, and note that $\tilde{Q}$ solves
\begin{eqnarray*}
\left\{ \begin{array}{c}
\left( \frac{\displaystyle \partial}{\displaystyle \partial t} + C(t) \right) \tilde{Q}(t) = -R(t), \\
\tilde{Q}(0) = \O (1)_{L^2 \to L^2}.
\end{array} \right.
\end{eqnarray*}
Since $R(t) = \O(1)_{L^2 \to L^2}$ by the induction hypothesis, Picard's theorem implies $\tilde{Q}(t):L^2 \to L^2$ 
as desired.
\end{proof}
\end{proof}

We will need to review some basic facts about the calculus of symbols with two parameters.  We will only use symbol spaces with two
parameters in the context of microlocal estimates, in which case we
may assume we are working in an open subset of $\reals^{2n}$.  We define the 
following spaces of symbols with two parameters:
\begin{eqnarray*}
\lefteqn{\s^{k,m, \widetilde{m}} \left( \reals^{2n} \right) := } \\
& = &\Big\{ a \in \Ci \left( \reals^{2n} \times (0,1]^2 \right):   \\ 
&& \quad \quad  \left| \partial_x^\alpha \partial_\xi^\beta a(x, \xi; h, \tilde{h}) \right| 
\leq C_{\alpha \beta}h^{-m}\tilde{h}^{-\widetilde{m}} \langle \xi \rangle^{k - |\beta|} \Big\}.
\end{eqnarray*}
For the applications in this paper, we assume $\tilde{h} >h$ and
define the scaled spaces:
\begin{eqnarray*}
\lefteqn{\s_{\delta}^{k,m, \widetilde{m}} \left(\reals^{2n} \right):= } \\
& = & \Bigg\{ a \in \Ci \left(\reals^{2n} \times (0,1]^2 \right):  \\
 && \quad \quad  \left| \partial_x^\alpha \partial_\xi^\beta a(x, \xi; h, \tilde{h}) \right| 
\leq C_{\alpha \beta}h^{-m}\tilde{h}^{-\widetilde{m}} \left( \frac{\tilde{h}}{h} \right)^{\delta(|\alpha| + |\beta|)} 
\langle \xi \rangle^{k - |\beta|} \Bigg\}.
\end{eqnarray*}
As before, we have the corresponding spaces of semiclassical pseudodifferential operators $\Psi^{k, m, \widetilde{m}}$ 
and $\Psi_{\delta}^{k,m, \widetilde{m}}$, where we will usually add a subscript of $h$ or $\tilde{h}$ to indicate which 
parameter is used in the quantization.  The relationship between $\Psi_h$ and 
$\Psi_{\tilde{h}}$ is given in the following lemma.  
\begin{lemma}
\label{U-hsc-lemma}
Let $a \in \s_{0}^{k,m,\tilde{m}}$, and set
\be
b(X, \Xi) = a\left(\csh^{\half}X, \csh^{\half} \Xi
\right) \in \s_{-\half}^{k,m,\tilde{m}}.
\ee
There is a linear operator $T_{h, \tilde{h}}$, unitary on $L^2$, and an operator 
such that
\be
\Op_{\tilde{h}}^w(b) T_{h, \tilde{h}} u = T_{h, \tilde{h}} \Op_h^w(a) u.
\ee
\end{lemma}
\begin{proof}
For $u \in L^2(\reals^{n})$, define $T_{h, \tilde{h}}$ by
\begin{eqnarray}
\label{U-hsc}
T_{h, \tilde{h}} u (X) := \csh^{\frac{n}{4}}u\left( \csh^{\half} X \right).
\end{eqnarray}
We see immediately that $T_{h, \tilde{h}}$ conjugates operators $a^w(x, hD_x)$ and $b^w(X, \tilde{h} D_X)$.
\end{proof}
We have the following microlocal commutator lemma.
\begin{lemma}
\label{2-param-lemma}
Suppose $a \in \s_{0}^{-\infty,0,0}$, $b \in
\s_{-\half}^{-\infty,m,\widetilde{m}}$, and $\tilde{h}>h$.  

(a) If $A = \Op_{\tilde{h}}^w (a)$ and $B = \Op_{\tilde{h}}^w(b)$,
\be
[A,B] &=& \frac{\tilde{h}}{i}\Op_{\tilde{h}}^w (\{a,b\}) + \O\left(
h^{3/2} \tilde{h}^{3/2} \right).
\ee

(b) More generally, for each $l>1$,
\be
 \ad_A^l B = \O_{L^2 \to L^2} \left( h \tilde{h}^{l-1} \right) .
\ee
\end{lemma}
\begin{proof}
Without loss of generality, $m = \widetilde{m} = 0$, so for (a) we have from the Weyl
calculus:
\be
[A,B] = \frac{\tilde{h}}{i} \Op_{\tilde{h}}^w (\{a,b\}) + \tilde{h}^3 \O
\left( \sum_{|\alpha| = |\beta| = 3} \partial^\alpha a \partial^\beta b \right),
\ee
since the second order term vanishes in the Weyl expansion of the
commutator.  Note $\partial^\alpha a$
is bounded for all $\alpha$, and observe for $|\beta|=3$,
\be
\tilde{h}^3 \partial^\beta b & = & \tilde{h}^3 \O \left(h^{3/2}
\tilde{h}^{-3/2} \right).
\ee
For part (b) we again assume $m = \widetilde{m} = 0$, and we observe that for $l>1$ we no longer have the
same gain in powers of $h$ as in part (a).  This follows from the fact
that the $\tilde{h}$-principal symbol for the commutator $[A,[A,B]]$, $-i
\tilde{h} \{a, -i \tilde{h}\{a,b\} \}$, satisfies 
\ben
-i \tilde{h}\{a, -i \tilde{h}\{a,b\}\} & = & -\tilde{h}^2 \Big( \partial_\Xi a
\partial_X \left(\partial_\Xi a \partial_X b - \partial_X a
\partial_\Xi b \right) \label{comm-calc-1} \\
&& \quad - \partial_X a \partial_\Xi \left( \partial_\Xi a \partial_X b - \partial_X a
\partial_\Xi b \right) \Big) \nonumber \\
& \in & \s_{0}^{-\infty,-1 ,-1}, \label{comm-calc-2}
\een
since $\{a,b\}$ involves products of derivatives of both $a$ and $b$.

For general $l>1$, assume 
\be
\sigma_{\tilde{h}} \left( \ad_A^l B \right) \in \s_{0}^{0, -1, 1-l}
\ee
and a calculation similar to (\ref{comm-calc-1}-\ref{comm-calc-2})
finishes the induction.
\end{proof}


\section{$h$-Fourier Integral Operators}
\label{FIO}
In this section we review some facts about $h$-Fourier Integral Operators ($h$-FIOs).  See \cite{duistermaat} for a 
comprehensive introduction to general FIOs without $h$, or
\cite{EvZw}, \S $10.1$ with the addition of the $h$ parameter.  
For this note, we are only interested in a special class of $h$-FIOs, namely those associated to a symplectomorphism.  
In order to motivate this, suppose $f:X \to Y$ is a diffeomorphism.  Then we write
\begin{eqnarray*}
f^*u(x) = u(f(x))= \frac{1}{(2 \pi h)^n} \int e^{i \langle f(x) -y, \xi \rangle /h} u(y) dy d\xi,
\end{eqnarray*}
and $f^*: \Ci(Y) \to \Ci(X)$ is an $h$-FIO associated to the nondegenerate phase function $\phi = \langle f(x) -y, 
\xi \rangle$.  We recall the notation from \cite{duistermaat}: if $A:\Ci_c(Y) \to \mathcal{D}'(X)$ is a continuous mapping with
distributional kernel $K_A \in \mathcal{D}'(X \times Y)$, 
\begin{eqnarray*}
\WF'(A) & = & \{((x, \xi),(y, \eta)) \in (T^*X \times T^*Y) \setminus 0 :
\\
&& \quad \quad (x, y ; \xi, - \eta) \in \WF (K_A) \}.
\end{eqnarray*}
In this notation, we note 
\begin{eqnarray*}
\WF' f^* \subset  \left\{ ((x, \xi),(y, \eta)): y = f(x), \,\, \xi = \,^t D_x f \cdot \eta \right\},
\end{eqnarray*}
which is the graph of the induced symplectomorphism
\begin{eqnarray*}
\kappa (x, \xi) = (f(x), (\,^tD_x f)^{-1}(\xi)).
\end{eqnarray*}
\indent To continue, we follow \cite{SjZw}, and let $A(t)$ be a smooth family of pseudodifferential operators: 
$A(t) = \Op_h^w(a(t))$ with 
\begin{eqnarray*}
a(t) \in \Ci \left( [-1,1]_t; \s^{-\infty,0} \left( T^*X \right) \right),
\end{eqnarray*}
such that for each $t$, $\WF(A(t)) \Subset T^*X$.  Let $U(t): L^2(X) \to L^2(X)$ be defined by
\begin{eqnarray}
\left\{ \begin{array}{c}
hD_tU(t) + U(t)A(t) = 0, \\
U(0) = U_0 \in \Psi_h^{0,0}(X),
\end{array} \right. \label{U(t)}
\end{eqnarray}
where $D_t = -i \partial / \partial t$ as usual.  If we let $a_0(t)$ be the real-valued $h$-principal symbol of 
$A(t)$ and let $\kappa(t)$ be the family of symplectomorphisms defined by
\begin{eqnarray*}
\left\{ \begin{array}{c}
\frac{\displaystyle d}{\displaystyle dt} \kappa(t)(x, \xi) = \left( \kappa(t) \right)_* \left( H_{a_0(t)}(x, \xi) \right), \\
\kappa(0)(x, \xi) = (x, \xi),
\end{array} \right.
\end{eqnarray*}
for $(x, \xi) \in T^*X$, then $U(t)$ is a family of $h$-FIOs associated to $\kappa(t)$.  We have the following 
well-known theorem of Egorov (see, for example \cite{EvZw}, \S $10.1$).
\begin{old-thm}[Egorov's Theorem]
Suppose $B \in \Psi_h^{k,m}(X)$, and $U(t)$ defined as above.  Suppose further that $U_0$ in \eqref{U(t)} is 
elliptic ($\sigma_h(U_0) \geq c >0$).  Then there exists a smooth family of pseudodifferential operators $V(t)$ 
such that
\begin{eqnarray}
\left\{ \begin{array}{c}
\sigma_h \left( V(t) B U(t) \right) = \left( \kappa(t) \right)^* \sigma_h(B), \\
V(t)U(t) -I, \,\, U(t)V(t) -I \in \Psi_h^{-\infty, -\infty}(X). \label{egorov1}
\end{array} \right.
\end{eqnarray}
\end{old-thm}
\begin{proof} As $U_0$ is elliptic, there exists an approximate inverse $V_0$, such that $U_0 V_0 -I, \,\, V_0 
U_0 - I \in \Psi_h^{-\infty, -\infty}$.  Let $V(t)$ solve
\begin{eqnarray*}
\left\{ \begin{array}{c}
h D_t V(t) - A(t) V(t) = 0, \\
V(0) = V_0.
\end{array} \right.
\end{eqnarray*}
Write $B(t) = V(t)BU(t)$, so that
\begin{eqnarray*}
hD_t B(t) = A(t)V(t)BU(t) - V(t)BU(t)A(t) = [A(t), B(t)] 
\end{eqnarray*}
modulo $\Psi_h^{-\infty, -\infty}$.  But the principal symbol of $[A(t), B(t)]$ is 
\begin{eqnarray*}
\sigma_h \left( [A(t), B(t)] \right) = \frac{h}{i} \left\{ \sigma_h ( A(t)), \sigma_h (B(t)) \right\} =
 \frac{h}{i} H_{a_0(t)} \sigma_h(B(t)),
\end{eqnarray*}
so \eqref{egorov1} follows from the definition of $\kappa(t)$.
\end{proof}
Let $U:= U(1)$, and suppose the graph of $\kappa$ is denoted by $C$.  Then we introduce the standard notation
\begin{eqnarray*}
U \in I_h^0(X \times X; C'), \,\,\, \text{with}\,\,\, C' = \left\{ (x, \xi; y, -\eta) : (x, \xi) = 
\kappa(y, \eta) \right\},
\end{eqnarray*}
meaning $U$ is the $h$-FIO associated to the graph of $\kappa$.  The next few results when taken together 
will say that locally all $h$-FIOs associated to symplectic graphs are of the same form as $U(1)$.  First
 a well-known lemma.
\begin{lemma}
\label{deform-lemma}
Suppose $\kappa : \nbhd (0,0) \to \nbhd (0,0)$ is a symplectomorphism fixing $(0,0)$.  Then there exists 
a smooth family of symplectomorphisms $\kappa_t$ fixing $(0,0)$ such that $\kappa_0 = \id$ and $\kappa_1 
= \kappa$.  Further, there is a smooth family of functions $g_t$ such that 
\begin{eqnarray*}
\frac{d}{dt} \kappa_t = (\kappa_t)_* H_{g_t}.
\end{eqnarray*}
\end{lemma}
The proof of Lemma \ref{deform-lemma} is standard, but we include a
sketch here, as it will be used in the proof of Proposition
\ref{normal-prop-1} (see \cite{EvZw} \S $10.1$ for details).
\begin{proof}[Sketch of Proof]
First suppose $K: \reals^{2n} \to \reals^{2n}$ is a linear symplectic
transformation.  Write the polar decomposition of $K$, $K = QP$ with
$Q$ orthogonal and $P$ positive definite.  It is standard that $K$
symplectic implies $Q$ and $P$ are both symplectic as well.  Identify
$\reals^{2n}$ with $\cx^n$ on which $Q$ is unitary.  Write $Q = \exp
iB$ for $B$ Hermitian and $P = \exp A$ for $A$ real symmetric and $JA
+ AJ = 0$, where
\begin{eqnarray*}
J := \left( \begin{array}{cc} 0 & -I \\ I & 0 \end{array} \right)
\end{eqnarray*}
is the standard matrix of symplectic structure on $\reals^{2n}$.
  Then $K_t = \exp (itB) \exp(tA)$ satisfies $K_0 = \id$ and
$K_1 = K$.

In the case $\kappa$ is nonlinear, set $K = \partial \kappa(0,0)$ and
choose $K_t$ such that $K_0 = \id$ and $K_\half = K$.  Then set 
\begin{eqnarray*}
\tilde{\kappa}_t(x, \xi) = \frac{1}{t} \kappa(t(x, \xi)),
\end{eqnarray*}
and note that $\tilde{\kappa}_t$ satisfies $\tilde{\kappa}_0 = K$,
$\tilde{\kappa}_1 = \kappa$.  Rescale $\tilde{\kappa}_t$ in $t$, so
that $\tilde{\kappa}_t \equiv K$ near $1/2$ and $\tilde{\kappa}_1 =
\kappa$.  Rescale $K_t$ so that $K_0 = \id$ and $K_t \equiv K$ near
$1/2$.  Then $\kappa_t$ is defined for $0 \leq t \leq 1$ by taking
$K_t$ for $0 \leq t \leq 1/2$ and $\tilde{\kappa}_t$ for $1/2 \leq t
\leq 1$.

To show $\frac{d}{dt} \kappa_t = (\kappa_t)_* H_{g_t}$, set $V_t =\frac{d}{dt}
\kappa_t$.  Cartan's formula then gives for $\omega$ the symplectic form
\begin{eqnarray*}
\mathcal{L}_{V_t} \omega = d\omega \contraction V_t + d( \omega \contraction
V_t),
\end{eqnarray*}
but $\mathcal{L}_{V_t} \omega = \frac{d}{dt}  \kappa_t^* \omega = 0$ since
$\kappa_t$ is symplectic for each $t$.  Hence $\omega \contraction V_t =
dg_t$ for some smooth function $g_t$ by the Poincar\'{e} lemma, in
other words, $V_t = (\kappa_t)_* H_{g_t}$.
\end{proof}

We have the following version of Egorov's theorem.
\begin{proposition}
\label{AF=FB}
Suppose $U$ is an open neighbourhood of $(0,0)$ and $\kappa: U \to U$ is a symplectomorphism fixing $(0,0)$.  
Then there is a bounded operator $F : L^2 \to L^2$ such that for all $A = \Op_h^w(a)$,
\begin{eqnarray*}
AF = FB \,\, \text{microlocally on}\,\, U \times U,
\end{eqnarray*}
where $B = \Op_h^w(b)$ for a Weyl symbol $b$ satisfying
\begin{eqnarray*}
b = \kappa^* a + \O(h^2).
\end{eqnarray*}
$F$ is microlocally invertible in $U \times U$ and $F^{-1} A F = B$ microlocally in $U \times U$.
\end{proposition}
Proposition \ref{AF=FB} is a standard result, however we include a proof as we will be using it for the proof 
of Theorem \ref{gamma-egorov}.
\begin{proof}
For $0 \leq t \leq 1$ let $\kappa_t$ be a smooth family of symplectomorphisms satisfying $\kappa_0 = \id$, 
$\kappa_1 = \kappa$, and let $g_t$ satisfy $\frac{d}{dt} \kappa_t =
(\kappa_t)_* H_{g_t}$.  Let $G_t = \Op_h^w(g_t)$, 
and solve the following equations 
\begin{eqnarray*}
&& \left\{ \begin{array}{c}
h D_tF(t) + F(t) G(t) = 0, \,\, (0 \leq t \leq 1) \\
F(0) = I,
\end{array} \right. \\
&& \left\{ \begin{array}{c}
hD_t \tilde{F}(t) - G(t) \tilde{F}(t) = 0,\,\, (0 \leq t \leq 1) \\
\tilde{F}(0) = I.
\end{array} \right.
\end{eqnarray*}
Then $F(t), \tilde{F}(t) = \O(1) : L^2 \to L^2$ and
\begin{eqnarray*}
hD_t\left( F(t) \tilde{F}(t) \right) = -F(t)G(t) \tilde{F}(t) + F(t) G(t) \tilde{F}(t) = 0,
\end{eqnarray*}
so $F(t) \tilde{F}(t) = I$ for $0 \leq t \leq 1$.  Similarly,
$E(t) = \tilde{F} F - I$ satisfies
\ben
\label{ode-eg}
hD_t E(t)  =  G(t) \tilde{F}(t) F(t) - \tilde{F}(t) F(t) G(t)  =  [G(t), E(t) ]
\een
with $E(0) = 0$.  But equation \eqref{ode-eg} has
unique solution $E(t) \equiv 0$ for the initial condition $E(0) = 0$.
Hence $\tilde{F}(t) F(t) =I$ microlocally. \\
\indent Now set $B(t) = \tilde{F}(t) A F(t)$.  We would like to show $B(t) = \Op_h^w(b_t)$, for $b_t = 
\kappa_t^*a + \O(h^2)$.  Set $\tilde{B}(t) = \Op_h^w (\kappa_t^*a)$.  Then
\begin{eqnarray*}
hD_t \tilde{B}(t) & = & \frac{h}{i} \Op_h^w \left( \frac{d}{dt}\kappa_t^*a \right) \\
& = & \frac{h}{i} \Op_h^w \left(\{ g_t, \kappa_t^*a \}\right) \\
& = & \left[ G(t), \tilde{B}(t) \right] + E_1(t),
\end{eqnarray*}
where $E_1(t) = \Op_h^w(e_1(t))$ for $e_1(t)$ a smooth family of symbols.  
Note if we take $g_t \# (\kappa_t^* a) - (\kappa_t^* a) \# g_t$, the composition formula \eqref{a-pound-b} 
implies the $h^2$ term vanishes for the Weyl calculus since $\omega^2$ is symmetric while 
\begin{eqnarray*}
g_t(x,\xi) \kappa_t^* a(y, \eta) - \kappa_t^* a(x, \xi) g_t(y, \eta)
\end{eqnarray*}
is antisymmetric.  Thus $E_1(t) \in \Psi_h^{0,-3}$, since we are working
microlocally.  We calculate
\begin{eqnarray}
\lefteqn{hD_t \left( F(t) \tilde{B}(t) \tilde{F}(t) \right) =} \label{B-de1}\\
& = & -F(t) G(t) \tilde{B}(t) \tilde{F}(t) + F(t) 
\left( \left[ G(t), \tilde{B}(t) \right] + E_1(t) \right) \tilde{F}(t) \label{B-de2}\\ 
&& \quad + F(t)\tilde{B}(t) G(t) \tilde{F}(t) \nonumber \\
& = & F(t) E_1(t) \tilde{F}(t) \label{B-de3}\\
& = & \O(h^3).\nonumber
\end{eqnarray}
Integrating in $t$ and dividing by $h$ we get
\begin{eqnarray}
\label{B-int}
F(t) \tilde{B}(t) \tilde{F}(t) = A + \frac{i}{h} \int_0^t F(s)E_1(s)\tilde{F}(s) ds = A + \O(h^2),
\end{eqnarray}
so that $\tilde{B}(t) - B(t) = \O(h^2)$. 

We will construct families of pseudodifferential operators $B_k(t)$ so
that for each $m$
\ben
B(t) = \tilde{B}(t) + B_1(t) + \cdots + B_m(t) + \O(h^{m+2}). \label{B-ind}
\een
Let 
\be
\tilde{e}_1(t) = (\kappa_t)^* \int_0^t (\kappa_s^{-1})^* e_1(s) ds,
\ee
and set $\tilde{E}_1(t) = \Op_h^w( \tilde{e}_1(t))$.  Observe
\be
hD_t \tilde{E}_1 = \left[ G(t), \tilde{E}_1 \right] + \frac{h}{i}\left(
E_1(t) + E_2(t)\right) ,
\ee
where $E_2(t) \in \Psi_h^{0,-4}$ by the Weyl calculus, since
$[G,\tilde{E}_1] = \O(h^4)$.  Then as in (\ref{B-de1}-\ref{B-de3})
\be
hD_t \left( F(t)\tilde{E}_1(t) \tilde{F}(t) \right) & = & - F(t)
\left[ G(t), \tilde{E}_1(t) \right] \tilde{F}(t) + F(t) hD_t\left(
\tilde{E}_1(t) \right) \tilde{F}(t) \\
& = & \frac{h}{i}\left( F(t) E_1(t) \tilde{F}(t) + F(t) E_2(t)
\tilde{F}(t)\right) .
\ee
Integrating in $t$ gives
\be
F(t) \tilde{E}_1(t) \tilde{F}(t) = \int_0^t F(s)E_1(s) \tilde{F}(s) ds
+ \frac{i}{h}\int_0^t F(s)E_2(s) \tilde{F}(s) ds,
\ee
and substituting in \eqref{B-int} gives
\be
\tilde{B}(t) - B(t) & = & \frac{i}{h} \tilde{E}_1(t) - \tilde{F}(t) \left(
\frac{i}{h} \int_0^t F(s)E_2(s)\tilde{F}(s)ds \right) F(t) \\
& = & \frac{i}{h} \tilde{E}_1(t) + \O(h^3).
\ee
Setting $B_1(t) = i \tilde{E}_1(t)/h$ and continuing inductively gives
$B_k(t)$ satisfying \eqref{B-ind}. 

\indent Let $l$ be a linear symbol, and $L = \Op_h^w(l)$.  Then
\begin{eqnarray*}
\ad_L (\tilde{B} - B) = \left[ \tilde{B} - B, L \right] = \O(h^2).
\end{eqnarray*}
Fix $N$.  From \eqref{B-ind} we can choose $B_1, \ldots , B_N$ so that
replacing $\tilde{B}$ with $\tilde{B}+ B_1 + \cdots + B_N$, we have for $l_1, \ldots, l_N$ linear symbols, $L_k = \Op_h^w(l_k)$,
\begin{eqnarray*}
\ad_{L_1} \circ \cdots \circ \ad_{L_N} (\tilde{B} - B) = \O(h^{N+2}),
\end{eqnarray*}
so Beals's Theorem implies $B(t) = \Op_h^w(b(t))$ for $b(t) = \kappa_t^* a + \O(h^2)$.
\end{proof}
The next proposition is essentially a converse to Proposition \ref{AF=FB}.
\begin{proposition}
\label{U-FIO}
Suppose $U = \O(1): L^2 \to L^2$ and for all pseudodifferential operators $A ,B \in \Psi_h^{0,0}(X)$ 
such that $\sigma_h(B) = \kappa^* \sigma_h(A)$, $AU = UB$ microlocally near $(\rho_0, \rho_0)$, 
where $\kappa : \nbhd (\rho_0, \rho_0) \to \nbhd (\rho_0 , \rho_0)$ is a symplectomorphism fixing 
$(\rho_0, \rho_0)$.  Then $U \in I_h^0(X \times X; C')$ microlocally near $(\rho_0, \rho_0)$.
\end{proposition}
\begin{proof}
Choose $\kappa_t$ a smooth family of symplectomorphisms such that $\kappa_0 = \id$, $\kappa_1 = \kappa$, 
and $\kappa_t( \rho_0) = \rho_0$.  Choose $a(t)$ a smooth family of functions satisfying $\frac{d}{dt} 
\kappa_t = (\kappa_t)_* H_{a(t)}$, and let $A(t) = \Op_h^w(a(t))$.  Let $U(t)$ be a solution to 
\begin{eqnarray*}
\left\{ \begin{array}{c} 
hD_t U(t) - U(t)A(t) = 0, \\
U(1) = U,
\end{array} \right.
\end{eqnarray*}
for $0 \leq t \leq 1$.  Next let $A$ and $B$ satisfy the assumptions of the proposition.  Since $AU = UB$, 
we can find $V(t)$ satisfying
\begin{eqnarray}
\label{V(t)-1}
\left\{ \begin{array}{c}
A U(t)V(t) = U(t)BV(t), \\
V(0) = \id.
\end{array} \right.
\end{eqnarray}
By Egorov's theorem, the right hand side of \eqref{V(t)-1} is equal to 
\begin{eqnarray*}
U(t) V(t) \left( V(t)^{-1} B V(t) \right) = U(t) V(t) A + \O(h).
\end{eqnarray*}
Setting $t = 0$, we see $[U(0), A ] = \O(h)$.  Applying the same argument to $[U(t), A]$ and another choice 
of $\tilde{A}, \tilde{B}$ satisfying the hypotheses of the proposition yields by induction,
\begin{eqnarray}
\label{beals-00}
\ad_{A_1} \circ \cdots \circ \ad_{A_N} U(0) = \O(h^N)
\end{eqnarray}
for any choice of $A_1, \ldots, A_N \in \Psi_h^{0,0}(X)$.  Since we are only interested in what $U(t)$ looks 
like microlocally, \eqref{beals-00} is sufficient to apply Beals's Theorem and conclude that $U(0) \in 
\Psi_h^{0,0}(X)$.  Thus $U(t)$ and hence $U(1) = U$ is in $I_h^0(X \times X; C')$ for the twisted graph
\begin{eqnarray*}
C' = \left\{ (x,\xi, y, -\eta): (y, \eta) = \kappa(x, \xi) \right\}.
\end{eqnarray*}
\end{proof}
Using the following more general version of the Poincar{\'e} lemma
from \cite{Wei2}, we will be able to generalize Proposition \ref{AF=FB} to a
neighbourhood of a periodic orbit.
\begin{lemma}
\label{poincare-lemma}
Let $N \subset T^*X$ be a closed submanifold, and assume $(x, \xi) \in
N$ implies $(x,0) \in N$.  Then if $\omega$ is a closed $k$-form such that $\left. \omega \right|_N = 0$, then there is a $(k-1)$-form $I( \omega)$ in a neighbourhood of $N$ such that $\omega = d I(\omega)$.
\end{lemma}
\begin{proof}
Let $m_s : T^*X \to T^*X$, $m_s: (x, \xi) \mapsto (x, s \xi)$, be multiplication by $s$ in the fibres, and define 
\begin{eqnarray*}
X_s = \left. \left( \frac{d}{dr} m_r \right) \right|_{r=s}.
\end{eqnarray*}
That is, in coordinates, 
\begin{eqnarray*}
X_s = \frac{1}{s} \sum_j \xi_j \frac{\partial}{\partial_{\xi_j}}
\end{eqnarray*}
is just $1/s$ times the radial vector field.  Then
\begin{eqnarray*}
\left. \frac{d}{dr} (m_r^* \omega) \right|_{r=s} = m_s^* \left( X_s \contraction d \omega \right) 
+ d \left( m_s^*( X \contraction \omega) \right),
\end{eqnarray*}
and integrating in $r$ gives
\begin{eqnarray*}
\omega - m_0^* \omega = I ( d\omega) + d I(\omega)
\end{eqnarray*}
for 
\begin{eqnarray*}
I (\omega) = \int_0^1 m_r^* (X_r \contraction  \omega ) dr.
\end{eqnarray*}
Now $\left. \omega \right|_N = 0$ and $ d \omega = 0$ finishes the proof.
\end{proof}

\begin{theorem}
\label{gamma-egorov}
Suppose $N \subset T^*X$ is a closed submanifold such that $(x, \xi)
\in N$ implies $(x,0) \in N$, and assume 
$\kappa: \neigh (N)  \to \kappa(\neigh(N))$ is a symplectomorphism which is smoothly homotopic in the symplectic 
group to identity on $N$.  Then there is a bounded linear operator $F: L^2(\neigh(N)) \to L^2(\kappa(\neigh(N)))$ 
such that for all $A = \Op_h^w(a)$,
\begin{eqnarray*}
AF = FB \,\,\, \text{microlocally on } \neigh(N) \times \kappa(\neigh(N)),
\end{eqnarray*}
where $B = \Op_h^w(b)$ for a Weyl symbol $b = \kappa^*a + \O(h^2)$.  Further, $F$ is microlocally invertible 
and $F^{-1}A F = B$ in $N \times \kappa(N)$.
\end{theorem}
\begin{proof}
The proof will follow from the proof of Proposition \ref{AF=FB}.  Let $\kappa_t$ be the homotopy in the 
Proposition, $\kappa_0 = \id$ and $\kappa_1 = \kappa$.  We need only verify that $\kappa_t$ is generated 
by a Hamiltonian.  Set $V_t = \frac{d}{dt} \kappa_t$, and calculate
\begin{eqnarray*}
0 = \frac{d}{dt} \kappa_t^* \omega = \mathcal{L}_{V_t} \omega = V_t \contraction d \omega + d ( V_t \contraction \omega ).
\end{eqnarray*}
Hence $\lambda_t = V_t \contraction \omega$ is closed and further $\left. \lambda_t \right|_{N} = 0$ so we may 
apply Lemma \ref{poincare-lemma} to obtain a $0$-form $I( \lambda_t)$ so that
\begin{eqnarray*}
d I(\lambda_t) = \lambda_t,
\end{eqnarray*}
or 
\begin{eqnarray*}
V_t = H_{I(\lambda_t)}.
\end{eqnarray*}
\end{proof}

We will make use of the following proposition (see \cite{EvZw} \S
$10.5$ for a proof).
\begin{proposition}
\label{hDx-prop}
Let $P \in \Psi_h^{k,0}(X)$ be a semiclassical operator of real principal type ($p = \sigma_h(P)$ is real 
and independent of $h$), and assume $dp \neq 0 $ whenever $p=0$.  Then for any $\rho_0 \in \{ p^{-1}(0) \}$, 
there exists a symplectomorphism $\kappa : T^*X \to T^* \reals^n$ defined from a neighbourhood of $\rho_0$ 
to a neighbourhood of $(0,0)$ and an $h$-FIO $T$ associated to its graph such that \\
\indent (i) $\kappa^* \xi_1 = p$, \\
\indent (ii) $TP = h D_{x_1} T$ microlocally near $(\rho_0; (0,0))$, \\
\indent (iii) $T^{-1}$ exists microlocally near $((0,0); \rho_0)$.
\end{proposition}


\section{Symplectic Geometry and Quadratic Forms}
\label{symplectic-geometry}
We now return to the setup of the introduction.  Let $P(h)$ satisfy all the assumptions from \S \ref{intro}.  
The main tool at our disposal is to use symplectomorphisms to transform the Weyl principal symbol into a 
different Weyl principal symbol which is in a more tractible form.  Then by Propositions \ref{AF=FB} and 
\ref{U-FIO}, any estimates we prove about the quantization of the transformed principal symbol will hold for 
the original operator modulo $\O (h^2)$.  
  
It is classical (see, for example \cite{AbMa}) that using our 
assumptions on $p$, the Implicit Function Theorem guarantees that there is an $\epsilon_0 >0$ such that for 
$\epsilon \in [-\epsilon_0, \epsilon_0]$, the energy surface $\{p^{-1}(\epsilon)\}$ is regular and contains a 
closed loxodromic orbit $\gamma^\epsilon$.  Further, 
\begin{eqnarray*}
\overline{\gamma} := \bigcup_{-\epsilon_0 \leq \epsilon \leq \epsilon_0} \gamma^\epsilon
\end{eqnarray*}
is a smooth, 2-dimensional symplectic manifold diffeomorphic to $\SS^1
\times [-\epsilon_0, \epsilon_0] \subset T^* \SS^1$.  Choose symplectic
coordinates $(t, \tau, x, \xi)$ in a neighbourhood of $\overline{\gamma}$ so that $\gamma$ is the image of the unit circle, $\SS^1 \ni t \mapsto 
\gamma(t)$, $t$ parametrizes $\gamma^\epsilon$, and $\gamma = \{t, 0;
0,0 \}$.  In \cite{AbMa} it is shown that $S=\{ t=0\}$ is a contact manifold with the contact form 
$\tom_{(x, \xi)} = i^* \omega$, where $i: S \hookrightarrow X$ is the inclusion.  
Then the Poincar\'e map preserves $p$ and $\tom$, modulo a term encompassing the period shift for 
$\epsilon \in [- \epsilon_0, \epsilon_0]$ different from zero and $(x,
\xi) \neq (0,0)$.  This
motivates our next change of variables.  Similar to \cite{Sjo3}, we
observe that $\tau$ depends only on the energy surface in which
$\gamma^\epsilon$ lies: $\tau = g( \epsilon )$.  $H_p$ is tangent to
the energy surface $\{ p^{-1} ( \epsilon) \}$ for each $\epsilon \in
[- \epsilon_0, \epsilon_0]$, so that
\begin{eqnarray*}
&& \partial_t p(t, \tau, x, 0) = \partial_t p(t, \tau, 0, \xi) = 0,
  \,\,\, \text{and} \\
&& \partial_x p(t, \tau, 0,0) = 0, \,\, \partial_\xi p(t, \tau, 0,0) =
  0,
\end{eqnarray*}
so that
\begin{eqnarray*}
p(t, \tau, 0,0) = f( \tau ) \,\,\, \text{and } p(t, 0, x, \xi ) = f(0)
+ \O_t( x^2 + \xi^2).
\end{eqnarray*}
Thus, there exists a smooth nonvanishing function $a( t, \tau, x,
\xi)$ defined in a neighbourhood of $\overline{\gamma}$ such that
\begin{eqnarray*}
a(t, \tau, x, \xi) p(t, \tau, x, \xi) = f(\tau) + \O_t(x^2 + \xi^2).
\end{eqnarray*}


Since the Hamiltonian vector field of $p$, $H_p$ is tangent to 
$\{p=0\}$, we can choose a Poincar\'{e} section contained in $\{p=0\}$, that is, 
a $2n-2$ dimensional submanifold $N$, transverse 
to $H_p$ on $\{p=0\}$ centered at $\gamma(0)$.  Let $S : N \to N$ be the Poincar\'{e} (first return) map near 
$\gamma(0)$.  Note that $\omega = dt \wedge d \tau + \tom_{(x, \xi)}$ is the symplectic form on $T^*X$ in our 
choice of coordinates, so $S$ preserves the $(2n-2)$ dimensional symplectic form $\tom$ on $N$.  Thus $S$ is 
a symplectic mapping $N \to N$, with $S(0) = 0$.  That $\gamma$ is loxodromic means none of the eigenvalues of 
$dS(0)$ lie on the unit circle.  In this section for simplicity we
consider only the case where all the eigenvalues are distinct, (the general case is handled in \S
\ref{non-distinct-section}).  
We think of $dS(0)$ as the linearization of $S$ near $0 \in N$, with $N$ 
identified with $T_0N$ near $0$.  \\
\indent We want to put 
$p$ into a normal form in a neighbourhood of $\gamma$.  Inspiration for this construction comes from
\cite{guillemin} and \cite{Sjo3}.  Let $q(\rho)$ be defined near $0 \in N$ and quadratic 
such that $dS(0) = \exp H_q$.  Let $\kappa_t$ be a smooth family of symplectomorphisms such that $\kappa_0 = \id$ 
while $\kappa_1 = S$.  Then from the proof of Lemma \ref{deform-lemma}
we can find $q_t(\rho)$ defined near $0 \in N$ so that 
\begin{eqnarray*}
q_t(\rho) = q(\rho) + f_t(\rho)
\end{eqnarray*}
with $f_t(\rho) = \O_t( |\rho|^3)$ and
\begin{eqnarray*}
\frac{d}{dt} \kappa_t = (\kappa_t)_* H_{q_t}.
\end{eqnarray*}
\begin{remark}
Here we see the first obstacle to extending these techniques to
include negative real eigenvalues: We want to write $dS(0) = \exp H_q$
for a real quadratic form $q$.  But this is impossible for some linear
symplectic transformations with negative eigenvalues as the example 
\begin{eqnarray*}
dS(0) = \left( \begin{array}{cc}
  -e^{2} & 0 \\ 0 & -e^{-2}
\end{array} \right)
\end{eqnarray*}
shows.  Here $dS(0)$ is symplectic, but cannot be written as $\exp
H_q$ with $q$ real.  Roughly, negative eigenvalues may be realized
only by deforming
a family of symplectomorphisms $\kappa_t$ {\it through an elliptic
  component}.
\end{remark}

Set $\tilde{p}(s, \sigma, \rho) = \sigma + q_s(\rho)$.  We will show $p$ and $\tilde{p}$ 
are equivalent under a symplectic change of coordinates on the set
$p^{-1}(0)$.  Then since both $p$ and $\tilde{p}$ have nonvanishing
differentials, we can write 
\begin{eqnarray}
\label{p=ptilde}
\kappa^* p = b(t, \tau, x, \xi) \tilde{p} 
\end{eqnarray}
for a smooth, positive function $b$ and a symplectomorphism $\kappa$.  
Indeed, we claim 
\begin{eqnarray*}
\exp (tH_{p})(s, \sigma, \rho) =  \left( s+t, \sigma_t(\rho,s,\sigma),  \kappa_{t+s} \circ \kappa_s^{-1}(\rho)\right)
\end{eqnarray*}
for some $\sigma_t(s, \sigma, \rho)$, giving \eqref{p=ptilde}.  To see this, set
\begin{eqnarray*}
\Phi_t(s, \rho) := \left( s+t,  \kappa_{t+s} \circ \kappa_s^{-1}(\rho) \right).
\end{eqnarray*} 
We need to check that $\left. \Phi_t \right|_{N \times \SS^1}$ is a $1$-parameter group.  We compute
\begin{eqnarray*}
\left. \Phi_{t_1 + t_2} \right|_{N \times \SS^1} (  s, \rho ) = \left( s + t_1 + t_2,  \kappa_{t_1 + t_2 + s} 
\circ \kappa_{s}^{-1}(\rho) \right).
\end{eqnarray*}
But we check
\begin{eqnarray*}
\lefteqn{ \left. \Phi_{t_1} \right|_{N \times \SS^1} \circ \Phi_{N \times \SS^1} ( s, \rho) =} \\
& = & \left. \Phi_{t_1} \right|_{N \times \SS^1} \left( s+ t_2,  \kappa_{t_2 + s} \circ \kappa_s^{-1}(\rho) \right) \\
& = & \left( s+t_1+t_2 ,  \kappa_{t_1 + t_2 + s} \circ \kappa_{t_2 + s}^{-1}( \kappa_{t_2+s} \circ 
\left(\kappa_{t_2 + s} \circ \kappa_s^{-1}(\rho) \right) \right),
\end{eqnarray*}
so the group law holds.  We need only verify that $p$ and $\tilde{p}$ have the same Poincar\'{e} map, so we check:
\begin{eqnarray*}
\left. \left( \frac{d}{dt} \left. \Phi_t \right|_{N \times \SS^1}( s, \rho ) \right) \right|_{t=0} = 
\left(1, H_{q_s}(\rho) \right),
\end{eqnarray*}
which is clear.  Note this construction depends only on the Poincar\'{e} map $S$ and is unique up to symplectomorphism. \\
\indent Next we want to examine what form the quadratic part $q(\rho)$ can take.  The fact that
 $S(0) = 0$ implies we can write 
\begin{eqnarray}
\label{q-1}
q ( \rho) = \frac{1}{2} \langle q''(0) \rho, \rho \rangle .
\end{eqnarray}
Now we define the Hamilton matrix $B$ by
\begin{eqnarray}
\label{b-matrix}
q (\rho) = \frac{1}{2} \tom(\rho, B \rho)
\end{eqnarray}
so that the symplectic transpose of $B$, $\,^{\tom}B$, is equal to
$-B$.  Note that $B$ is the matrix representation of $H_q$, and so has
eigenvalues which are 
the logarithms (with a suitably chosen branch cut) of the eigenvalues
of $dS(0)$.  Thus the condition
 that $\gamma$ be loxodromic implies none of the eigenvalues of $B$ have nonzero real part.  Recall that 
since $dS(0)$ is a symplectic transformation, if $\mu$ is an eigenvalue of $dS(0)$, then so are 
$\overline{\mu}$, $\mu^{-1}$, and $\overline{\mu}^{-1}$.  This implies
for the corresponding Hamilton matrix $B$ 
in (\ref{b-matrix}), if $\lambda$ is an eigenvalue of $B$, then so are $-\lambda$, $\overline{\lambda}$, 
and $-\overline{\lambda}$.  Thus the analysis of $B$ in the loxodromic, or complex hyperbolic case amounts 
to analyzing the eigenvalues in sets of $2$ or $4$.  For this we
follow the appendix in \cite{IaSj}, and recall for this section we are
assuming the eigenvalues are distinct.  There are $2$ cases.  First, assume $\lambda_j > 0$ is real.  Then $-\lambda_j$ is also an eigenvalue.  
Let $e_j$ and $f_j$ be the respective eigenvectors such that $\widetilde{\omega}(e_j, f_j)=1$.  Then $e_j$ 
and $f_j$ span a real symplectic vector space of dimension $2$.  For a point $\rho$ in this vector space, 
write $\rho = x_je_j + \xi_j f_j$.  Then $(x_j, \xi_j)$ are symplectic coordinates, in which $q_j(\rho)$, 
the projection of $q$ onto the $j$th coordinates becomes $q_j( \rho) = \lambda_j x_j \xi_j$.  We call the
\begin{eqnarray*}
\lambda_j x_j \xi_j
\end{eqnarray*}
the {\it action variables}. \\  
\indent Now we would like to see what these actions look like when the eigenvalues have nonzero imaginary part.  
Suppose $\lambda_j$ is an eigenvalue with $\Re \lambda_j > 0$, $\Im \lambda_j >0$.  Then $-\lambda_j$, 
$\overline{\lambda}_j$, and $-\overline{\lambda}_j$ are eigenvalues.  Let $e_j$, $f_j$, $\overline{e}_j$, 
and $\overline{f}_j$ be the respective eigenvectors.  Note $\tom (e_j, \overline{e}_j)= \tom(e_j, \overline{f}_j) = 
\tom (f_j, \overline{f}_j) = 0$.  Scale $f_j$ so that $\tom (e_j, f_j) = 1$.  Then $\{e_j, f_j \}$ and $
\{\overline{e}_j, \overline{f}_j \}$ span complex conjugate symplectic vector spaces of complex dimension $2$.  
Thus $\{ e_j, \overline{e}_j, f_j, \overline{f}_j \}$ span a symplectic vector space of complex dimension $4$ 
which is the complexification of a real symplectic vector space.  Write a point $\rho$ in this space in this basis, 
$\rho = z_j e_j + \zeta_j f_j + w_j \overline{e}_j + \eta_j \overline{f}_j$.  Then $(z_j, \zeta_j, w_j, \eta_j)$ 
become symplectic coordinates, in which the projection $q_j$ becomes $q_j(\rho) = \lambda_j z_j \zeta_j + 
\overline{\lambda}_j w_j \eta_j$.  Now write 
\begin{eqnarray*}
e_j = \frac{1}{\sqrt{2}} \left( e_j^1 + i e_j^2 \right), \quad f_j = \frac{1}{\sqrt{2}} \left( f_j^1 - if_j^2 \right),
\end{eqnarray*}
for real $e_j^k$, $f_j^k$.  This is a symplectic change of basis, and writing $\rho$ in this basis:
\begin{eqnarray*}
\rho = z_j e_j + \zeta_j f_j + w_j \overline{e}_j + \eta_j \overline{f}_j = \sum_{k=1}^2 
\left( x_j^k e_j^k + \xi_j^k f_j^k \right),
\end{eqnarray*}
we have
\begin{eqnarray*}
q_j(\rho) = \Re \lambda_j \left( x_j^1 \xi_j^1 + x_j^2 \xi_j^2 \right) - \Im \lambda_j 
\left( x_j^1 \xi_j^2 - x_j^2 \xi_j^1 \right).
\end{eqnarray*}
This is summarized in the following proposition (using the notation of
\cite{IaSj}).  Let $n_{hc}$ be the number of complex hyperbolic eigenvalues
$\mu_j$ of $dS(0)$ with $|\mu_j|>1$, and $n_{hr}$ the number of real
hyperbolic eigenvalues $\mu_j$ of $dS(0)$ such that $\mu_j >1$.  Thus
we have $2n-2 = 4n_{hc} + 2 n_{hr}$.
\begin{proposition}
\label{normal-prop-1}
Let $p \in \Ci ( T^*X)$, $\gamma \subset \{p = 0 \}$ as in the introduction, with
the linearized Poincar\'{e} map having distinct eigenvalues $\mu_j$ not on the
unit circle.  Assume for $1 \leq j \leq n_{hc}$ we have $|\mu_j|>1$
and $\Im \mu_j >0$, and for $2n_{hc}+1 \leq j \leq 2n_{hc} + n_{hr}$
we have $\mu_j >1$.  Then there exists a neighbourhood, $U$, of $\gamma$ in $T^*X$, 
a smooth positive function $b \geq C^{-1} >0$ defined in $U$, and a 
symplectomorphism $\kappa : U \to \kappa(U) \subset T^*\SS_{(t, \tau)}^1 \times T^* \reals_{(x, \xi)}^{n-1}$ such that 
\begin{eqnarray*}
\kappa(\gamma)   =  \{ (t,0;0,0) : t \in \SS^1 \}, 
\end{eqnarray*}
and $b(t, \tau, x, \xi) p = \kappa^*(g + r)$, with
\begin{eqnarray}
\lefteqn{ g(t, \tau; x, \xi) =} \nonumber \\ 
& = &  \tau + \sum_{j=1}^{n_{hc}}\left( \Re \lambda_j \left( x_{2j-1} \xi_{2j-1} + x_{2j} \xi_{2j} \right) 
- \Im \lambda_j \left( x_{2j-1} \xi_{2j} - x_{2j}\xi_{2j-1} \right) \right) \label{p-pullback-def-1} \\
 & & \quad +  \sum_{j = 2n_{hc}+1}^{2n_{hc} + n_{hr}} \lambda_j x_j \xi_j , \quad \mathrm{with} \,\,\, 
2n_{hc} + n_{hr} = n-1 \,\,\, \mathrm{and} \label{p-pullback-def-2} \\
r & = & \O (|x|^3 + |\xi|^3). \nonumber
\end{eqnarray}
Here $\lambda_j = \log (\mu_j)$ for $|\mu_j|>1$ and $\Im \lambda_j
\geq 0$.
\end{proposition}
\begin{remark}
The quadratic form (\ref{p-pullback-def-1}-\ref{p-pullback-def-2}) in Proposition \ref{normal-prop-1} 
is the leading part of the real Birkhoff normal form for a symplectomorphism near a loxodromic fixed point.  
With a non-resonance condition and the addition of some higher order
``action'' variables (see, for example, \cite{hz} and \cite{IaSj}), the error $r$ could be taken to be 
\begin{eqnarray*}
r & = & \O (|x|^4 + |\xi|^4),
\end{eqnarray*}
or even $\O(|x|^\infty + |\xi|^\infty)$.
\end{remark}
\begin{remark}
We think of $p(t, \tau, x, \xi) \in \Ci( \reals^4)$, $p = \tau +
\lambda x \xi$, for $\lambda>0$ as our ``model case''.  The feature we
are going to exploit about this model case is that if $G(t, \tau, x,
\xi) = \half ( x^2 - \xi^2)$, then
\begin{eqnarray}
\label{hp-quad}
H_pG = \lambda ( x^2 + \xi^2),
\end{eqnarray} 
which is a positive definite quadratic form.  However,
the growth of $x^2 - \xi^2$ will force us to use instead $G(x, \xi) =
\log (1 + x^2) - \log(1 + \xi^2)$.  Suppose $p = \tau + \lambda x
\xi + x^3 -\xi^3 = \tau + \lambda x\xi + \O(x^3 + \xi^3)$ in a
neighbourhood of $\gamma$ of size $\epsilon >0$ as in
Proposition \ref{normal-prop-1}.  Then
\begin{eqnarray*}
H_p G = \lambda \frac{x^2}{1+ x^2} + \lambda \frac{\xi^2}{1 + \xi^2} +
3 \frac{\xi^2 x}{1 + x^2} + 3 \frac{ x^2 \xi}{1 + \xi^2}.
\end{eqnarray*}
Motivated by \eqref{hp-quad}, we would like to write this as
\begin{eqnarray*}
H_pG = \lambda \frac{x^2}{1+ x^2}(1 + \O( \epsilon)) + \lambda
\frac{\xi^2}{1 + \xi^2}(1 + \O(\epsilon)),
\end{eqnarray*}
which we clearly cannot do in this example.
\end{remark}

As the last remark indicates, in order to deal with the error terms, we will need a more refined form than that given in 
Proposition \ref{normal-prop-1}.  Inspiration for this development, and in particular Proposition \ref{normal-prop-2} 
comes from \cite{GeSj} and \cite{sjostrand}.  

Let $\{\mu_j\}$ be the eigenvalues of the linearized Poincar\'{e} map at $\gamma(0)$.  They come in pairs $\mu_j$, 
$\mu_j^{-1}$ for the real $\mu_j$ and in sets of four $\mu_j$, $\overline{\mu_j}$, $\mu_j^{-1}$, and 
$\overline{\mu}_j^{-1}$ for the complex $\mu_j$.  The Stable/Unstable Manifold Theorem guarantees we will get two 
$n$-dimensional, transversal, flow-invariant sub-manifolds $\Lambda_+$ and $\Lambda_-$ such that $\exp tH_p$ is 
expanding on $\Lambda_+$ and contracting on $\Lambda_-$.  Since the $\Lambda_\pm$ are invariant under the flow 
$\Phi_t = \exp tH_p$ which is symplectic, the symplectic form $\omega$ vanishes on the $\Lambda_\pm$, that is, 
the $\Lambda_\pm$ are Lagrangian submanifolds.  
\begin{lemma}
\label{zero-cov}
Assume $p$ is in the form of Proposition \ref{normal-prop-1}.  Then there exists a local symplectic coordinate system $(t,\tau,x, \xi)$ near $\gamma$ such that 
$\Lambda_+ = \{ \tau = 0, \xi = 0 \}$ and $\Lambda_- = \{ \tau = 0, x = 0 \}$.   
\end{lemma}
\begin{proof}
We claim the $\Lambda_\pm$ are orientable and embedded 
in $T^* \SS^1 \times T^* \reals^{n-1}$.  Since $dS(0)$ describes how the flow of $H_p$ has acted at time $t = 1$, 
we know the evolution of a tangent frame of $\Lambda_\pm$ will be
described by $dS(0)$.  Using the action variables in Proposition
\ref{normal-prop-1}, we have
\begin{eqnarray*}
dS(0) = \left( \begin{array}{cc}
  A & 0 \\0 & B \end{array} \right)
\end{eqnarray*}
with
\begin{eqnarray*}
A = \diag(  \mu_{1},
\bar{\mu}_{1}, \ldots \mu_{n_{hc}}, \bar{\mu}_{n_{hc}} ; \mu_{2n_{hc}+1}, \ldots, \mu_{2n_{hc}+n_{hr}}),
\end{eqnarray*}
describing the time $1$ evolution of $\Lambda_+$ and $| \mu_j | >1$ for each $1
\leq j \leq n_{hr}+n_{hc}$ by our choice of coordinates.  Similarly, 
\begin{eqnarray*}
B = \diag( \mu_{1}^{-1},
\bar{\mu}_{1}^{-1}, \ldots \mu_{n_{hc}}^{-1}, \bar{\mu}_{n_{hc}}^{-1} ; \mu_{2n_{hc}+1}^{-1}, \ldots, \mu_{2n_{hc}+n_{hr}}^{-1} )
\end{eqnarray*}
describes the time $1$ evolution of $\Lambda_-$ with $|\mu_j^{-1}| <
1$ for each $j$.  But we've assumed there are no negative real eigenvalues, so $\det A >0$ implies $\Lambda_+$ is
orientable.  Similarly, $\det B >0$ and $\Lambda_-$ is orientable.  Now our assumptions on $p$ mean the flow has no critical 
points in a neighbourhood of $\gamma$ so the $\Lambda_\pm$ can have no self intersections and hence are embedded.  \\
\indent Let $\tilde{\Lambda} \subset T^*\SS^1 \times T^* \reals^{n-1}$, $\tilde{\Lambda} = \{ \tau = 0, \xi=0 \}$.  
Since $\Lambda_+$ is a closed, $n$-dimensional submanifold of $T^*X$, the tubular neighbourhood theorem guarantees 
there is a diffeomorphism $f$ (not necessarily symplectic) taking a neighbourhood $U$ of $\gamma$ into itself so that 
$f$ fixes $t$ and
\begin{eqnarray*}
f ( \Lambda_+ \cap U ) = \tilde{\Lambda} \cap U.
\end{eqnarray*}
Further, since $T_{\gamma(t)}\Lambda_+ = T_{\gamma(t)}\tilde{\Lambda}$ for $0 \leq t \leq 1$, we can choose $f$ satisfying 
\begin{eqnarray}
\label{f-pb-nondeg}
\left[ (f^{-1})^* \widetilde{\omega} \right]_{\gamma(t)}= \widetilde{\omega}_{\gamma(t)}, \,\,\, 0 \leq t \leq 1.
\end{eqnarray}
The statement in the lemma about $\Lambda_+$ now follows directly from the more general Theorem 4.1 in \cite{Wei2}, but 
we include a proof of this concrete case.  We have $\tilde{\Lambda} \subset T^*\SS^1 \times T^* \reals^{n-1}$, a 
Lagrangian submanifold with two distinct symplectic structures, $\omega_0 = (f^{-1})^* \widetilde{\omega}$ and the 
standard symplectic structure $\omega_1$ inherited from $T^*\SS^1 \times T^* \reals^{n-1}$.  We want to find a 
diffeomorphism $g: U \to U$ such that $g( \tilde{\Lambda}) = \tilde{\Lambda}$ and $g^* \omega_1 = \omega_0$.  \\
\indent Set $\omega_s = s \omega_0 + (1 - s) \omega_1$.  We have $d \omega_s = 0$ and 
$\left. \omega_s \right|_{\tilde{\Lambda}} = 0$.  Note \eqref{f-pb-nondeg} implies $\omega_s$ is nondegenerate
 in a neighbourhood of $\gamma$ for $0 \leq s \leq 1$.  Let $\widehat{\omega}_s: TX \to T^*X$ denote 
the isomorphism generated by $\omega_s$, $\widehat{\omega}_s : Z \mapsto Z \contraction \omega_s$.  We use the general 
Poincar\'{e} Lemma \ref{poincare-lemma} to obtain a $1$-form $\phi = I(\omega_0 - \omega_1)$ so that 
$d\phi = \omega_0 - \omega_1$ and set $Y_s = \widehat{\omega}_s^{-1}( \phi )$.  Then 
$\left. \phi \right|_{\tilde{\Lambda}} = 0$ implies 
\begin{eqnarray*}
Y_s \contraction \omega_s & = & \widehat{\omega} ( Y_s ) \\
& = & \phi,
\end{eqnarray*}
so that $Y_s$ is tangent to $\tilde{\Lambda}$.  Thus if $g_s = \exp (s Y_s)$ for $0 \leq s \leq 1$ is the 
integral of $Y_s$, $g_s(\tilde{\Lambda}) = \tilde{\Lambda}$.  We calculate:
\begin{eqnarray*}
\frac{d}{dr} \left. \left( g_r^* \omega_r \right) \right|_{r = s} & = & g_s^* \left. 
\left( \frac{d}{dr} \omega_r \right) \right|_{r = s} + g_s^* \left( d ( Y_s \contraction \omega_s)\right) \\
& = & g_s^* \left( \omega_0 - \omega_1 + d( -\phi) \right) \\
& = & 0.
\end{eqnarray*}
Setting $g = g_1$ gives $g^* \omega_1 = \omega_0$ as desired.  Now taking $g^{-1} \circ f$ gives a 
diffeomorphism of a neighbourhood of $\gamma$ taking $\Lambda_+$ to $\tilde{\Lambda}$ such that 
$g^* \circ (f^{-1})^* \widetilde{\omega} = \widetilde{\omega}$. \\
\indent After this change of coordinates, we still need to put $\Lambda_-$ in the desired form.  
Since $\Lambda_-$ is transversal to $\Lambda_+$ and all of our transformations so far leave 
$\{ \tau = 0\}$ invariant, we can write $\Lambda_-$ as a graph over $\{x=0\}$:
\begin{eqnarray}
\label{lambda-minus-1}
\Lambda_- = \left\{ ( t, 0, x, \xi): x = g(\xi, t) \right\}.
\end{eqnarray}
Further, since for each fixed $t$, \eqref{lambda-minus-1} is Lagrangian and the first de Rham cohomology 
group $H^1_{dR}(\{\tau = 0,x=0\}) \simeq H^1_{dR}( \reals^{n-1})$ vanishes, it is classical that we can 
write $g(\xi, t) = \partial_\xi h(\xi, t)$ for a smooth $h(\xi, t)$ (see, for example, \cite{Lee}).  Then we write
\begin{eqnarray*}
\Lambda_- = \left\{ ( t, 0, x, \xi): x = \partial_\xi h(\xi, t) \right\},
\end{eqnarray*}
and observe $h$ must satisfy $\partial_\xi h(0,t) = 0$.  This determines $h$ up to a constant, which we take 
to be $0$ so that $h(0,t) = 0$.  Now let $b(\xi, t)$ be a smooth function satisfying $b(\xi, t) = 
\partial_t h(\xi, t)$, and note $b(0,t) = 0$.  Then we perform the following change of variables:
\begin{eqnarray*}
\left\{ \begin{array}{rcl} 
t' & = & t \\
\tau' & =& \tau + b(\xi,t) \\
x' & = & x - \partial_\xi h(\xi, t) \\
\xi' & = & \xi . \end{array} \right. 
\end{eqnarray*}
We calculate:
\begin{eqnarray*}
d \tau' \wedge dt' + d \xi' \wedge dx' & = & \left(d \tau + \sum_j \partial_{\xi_j} b(\xi, t) d\xi_j + 
\partial_t b(\xi, t) dt \right) \wedge dt \\
&&  + \sum_j d \xi_j \wedge \left(dx_j - \sum_i \partial_{\xi_i} \partial_{\xi_j} h(\xi ,t) d\xi_i - 
\partial_t \partial_{\xi_j} h(\xi, t) dt\right) \\
& = & d \tau \wedge d t + d \xi \wedge dx,
\end{eqnarray*}
by the symmetry of the Hessian $\partial_{\xi_i} \partial_{\xi_j} h(\xi, t)$.  Thus this change of 
variables is symplectic and the Lemma is proved.
\end{proof}

Using the change of variables in Lemma \ref{zero-cov}, we have the following proposition.  
\begin{proposition}
\label{normal-prop-2}
Let $p \in \Ci ( T^*X)$, $\gamma \subset \{p = 0 \}$ as above, with
the Poincar\'{e} map having distinct eigenvalues $\mu_j$ not on the
unit circle.  
Then there exists a neighbourhood, $U$, of $\gamma$ in $T^*X$, a smooth
positive function $b\geq C^{-1} >0$ defined in $U$, 
a symplectomorphism 
$\kappa : U \to \kappa(U) \subset T^*\SS_{(t, \tau)}^1 \times T^* \reals_{(x, \xi)}^{n-1}$, 
and a smooth, $n \times n$-matrix valued function $B_t$ such that 
\begin{eqnarray}
\kappa(\gamma)  & = & \{ (t,0;0,0) : t \in \SS^1 \}, \quad
\mathrm{and} \,\,\, b(t, \tau, x, \xi) p = \kappa^*g, \,\,\, 
\mathrm{with} \nonumber \\
g(t, \tau; x, \xi) & = &\tau + \langle B_t(x, \xi) x, \xi \rangle \label{normal-form-2},
\end{eqnarray}
with $B_t$ satisfying 
\begin{eqnarray}
\lefteqn{\left\langle B_t(0,0)x, \xi \right\rangle = } \nonumber \\ 
& = &  \sum_{j=1}^{n_{hc}}\left( \Re \lambda_j \left( x_{2j-1} \xi_{2j-1} + x_{2j} \xi_{2j} \right) 
- \Im \lambda_j \left( x_{2j-1} \xi_{2j} - x_{2j}\xi_{2j-1} \right) \right) \label{B-condition.a} \\
 & & +  \sum_{j = 2n_{hc}+1}^{2n_{hc} + n_{hr}} \lambda_j x_j \xi_j \label{B-condition.b}.
\end{eqnarray}
Here $\lambda_j = \log (\mu_j)$ for $|\mu_j|>1$ and $\Im \lambda_j
\geq 0$.
\end{proposition}
\begin{proof}
Recall that the Poincar\'{e} map $S$ is linear in lowest order, and let $dS(0)$ be the linearized map.  
Let $q_0$ satisfy $dS(0) = \exp H_{q_0}$.  After a linear symplectic change of variables, $q_0$ can be 
written in block-diagonal form 
\begin{eqnarray*}
q_0(x, \xi) & =&  \langle bx ,\xi\rangle \\
& = & \sum_{j=1}^{n_{hc}}\left( \Re \lambda_j \left( x_{2j-1} \xi_{2j-1} + x_{2j} \xi_{2j} \right) 
- \Im \lambda_j \left( x_{2j-1} \xi_{2j} - x_{2j}\xi_{2j-1} \right) \right)  \\
 & & +  \sum_{j = 2n_{hc}+1}^{2n_{hc} + n_{hr}} \lambda_j x_j \xi_j , \quad \mathrm{with} \,\,\, 2n_{hc} + n_{hr} = 2n-2.
\end{eqnarray*}
According to Lemma \ref{zero-cov}, we may symplectically change
variables so  $\Lambda_+ = \{\tau=0, \xi=0 \}$ and $\Lambda_- = \{ \tau=0, x=0\}$.  
The linearization of the Hamiltonian vector field of $p$ is $H_{q_0}$, which implies we have a quadratic 
form as in the proposition.
\end{proof}


\section{Proof of Theorem \ref{main-theorem-1} }
\label{main-theorem-1-proof}
\numberwithin{equation}{section}

\begin{proof}[Proof of Theorem \ref{main-theorem-1} with Distinct Eigenvalues]
First we assume $P(h)$ has principal symbol given by 
\begin{eqnarray}
\label{assume-p}
p(t, \tau; x, \xi) = \tau +\langle B_t(x, \xi) x, \xi \rangle,
\end{eqnarray}
with $B_t$ satisfying (\ref{B-condition.a}-\ref{B-condition.b}) as in Proposition \ref{normal-prop-2}.  Let $U$ be 
a neighbourhood of $\gamma$, $U \subset T^* \SS^1 \times T^* \reals^{n-1}$, and assume 
\begin{eqnarray*}
U \subset U_{\epsilon/2} := \left\{ (t, \tau, x, \xi): \left( d(x, x(\gamma(t)))^2 + |\xi - \xi(\gamma(t))|^2 
+ \tau^2 \right)^\half < \frac{\epsilon}{2} \right\}
\end{eqnarray*}
for $\epsilon >0$.  Let $\psi_0$ be a microlocal cutoff function to a neighbourhood of $U$, that is, take 
$\psi_0 \in \Ci_c(\reals^{2n})$, $\psi_0 \equiv 1$ on $U_{\epsilon/2}$ with support in $U_{ \epsilon}$.  
Then we assume throughout that we are working in $U_\epsilon$.  With $\tilde{h}$ small (fixed later 
in the proof), we do the following rescaling:
\begin{eqnarray}
\label{rescaling1}
X : = \hsc^\half x, \quad \Xi = \hsc^\half \xi.
\end{eqnarray}
and assume for the remainder of the proof that $|(X, \Xi)| \leq \hsc^\half \epsilon$.  We use the unitary operator 
$T_{h, \tilde{h}}$ defined in \eqref{U-hsc} to introduce the second parameter into $P(h)$.  Following \cite{BZ} we define the
 operator $\widetilde{P}(h)$ by 
\begin{eqnarray*}
\widetilde{P}(h) = T_{h, \tilde{h}} P(h) T_{h, \tilde{h}}^{-1},
\end{eqnarray*}
\noindent so that the principal symbol of $\widetilde{P}(h)$ is 
\begin{eqnarray}
\label{tildep-expression}
\lefteqn{\widetilde{p}(t, \tau; X, \Xi) = } \\
&& = \tau + \left\langle B_t\left(
\csh^{\half}( X, \Xi) \right) \csh^\half X, \csh^\half \Xi \right\rangle, \nonumber
\end{eqnarray}
\noindent and $\widetilde{p} \in \s_{-\half}^{-\infty,0,0}$ microlocally.  We have 
\begin{eqnarray}
\label{tildep-est}
\left|\partial_{X, \Xi}^\alpha \widetilde{p}\right| \leq C_\alpha\csh^{|\alpha|/2}
\end{eqnarray}
for $(X, \Xi) \in U_{\hsc^\half \epsilon}$ by Lemma \ref{U-hsc-lemma}. \\
\indent We will use the following escape function, which we define in the $(X, \Xi)$ coordinates:
\begin{eqnarray*}
G(X, \Xi) : = \frac{1}{2} \left( \log ( 1 + |X|^2) - \log( 1 + |\Xi |^2 ) \right).
\end{eqnarray*}
$G$ satisfies
\begin{eqnarray*}
\left| \partial_X^\alpha \partial_\Xi^\beta G(X, \Xi) \right| \leq C_{\alpha \beta} \langle X \rangle^{-|\alpha|}\langle 
\Xi \rangle^{-|\beta|}, \,\,\,\, \,\,\, \mathrm{for} \,\,\, (\alpha, \beta) \neq (0,0),
\end{eqnarray*}
and since $\langle X \rangle^2 \langle \Xi \rangle^{-2} $ is an order function, $G$ satisfies the assumptions of 
Lemma \ref{etG-lemma} so we may construct the family $e^{sG^w}$ for sufficiently small $s$.  \\
\indent Now for $|(X, \Xi)| \leq \hsc^{\half} \epsilon$ we have 
\begin{eqnarray}
\lefteqn{ H_{\widetilde{p}} G(X, \Xi) = } \nonumber \\
& = & \csh\left[ \left\langle B_t X, \frac{\partial}{\partial X} \right\rangle
  - \left\langle B_t \frac{\partial}{\partial \Xi}, \Xi \right\rangle \right]  G(X, \Xi) \label{hpg-main-term} \\
& & \quad + \csh^{\frac{3}{2}}\left[ \sum_{j = 1}^{n-1} \left\langle \frac{\partial}{\partial \Xi_j}B_t( \cdot, \cdot ) 
X, \Xi \right\rangle \frac{\partial}{\partial X_j}G(X, \Xi)\right] \label{hpg-err1} \\
& & \quad - \csh^{\frac{3}{2}}\left[ \sum_{j = 1}^{n-1} \left\langle \frac{\partial}{\partial X_j} B_t( \cdot,\cdot ) X, 
\Xi \right\rangle \frac{\partial}{\partial \Xi_j}  G(X, \Xi)\right]. \label{hpg-err2}
\end{eqnarray}
\indent For $s$ sufficiently small, we define a family of operators
\begin{eqnarray}
\label{Pth1}
\widetilde{P}_s(h) &=& e^{-s G^w} \widetilde{P}(h) \Op_{\tilde{h}}^w \left(\psi_0\left(\csh^{\half} \bullet \right) \right) 
e^{sG^w} \nonumber \\
& = & \exp \left(-s \ad_{G^w}\right) \widetilde{P}(h) \Op_{\tilde{h}}^w \left( \psi_0 \left(\csh^{\half} \bullet \right) \right),
\end{eqnarray}
\noindent where $\Op_{\tilde{h}}^w$ and $G^w$ are quantizations in the $\tilde{h}$-Weyl calculus.  Now owing to 
Lemma \ref{2-param-lemma} and
\eqref{tildep-est} we have microlocally to leading order in $h$:
\begin{eqnarray*}
\ad_{G^w}^k \left( \widetilde{P}(h)
\Op_{\tilde{h}}^w \left(\psi_0\left(\csh^{\half} \bullet
\right)\right) \right) = \O_{L^2 \to L^2}\left( h \tilde{h}^{k-1} \right),
\end{eqnarray*}
and in particular,
\begin{eqnarray}
\label{PG-comm}
\left[ \widetilde{P}(h), G^w \right]  =  -i \tilde{h} \Op_{\tilde{h}}^w \left( H_{\widetilde{p}} G\right) 
+ \O (h^{3/2} \tilde{h}^{3/2}).
\end{eqnarray}
\indent Now near $(0,0)$, $B_t$ is positive definite, $\langle B_t X,X \rangle \geq C^{-1} |X|^2$, so 
\begin{eqnarray*}
\left\langle B_t X,X \right\rangle^{-1} \leq C|X|^{-2}.
\end{eqnarray*}
Applying this to the errors (\ref{hpg-err1}-\ref{hpg-err2}) we get
\begin{eqnarray*}
\csh^{\frac{3}{2}}\left[ \sum_{j = 1}^{n-1} \left\langle \frac{\partial}{\partial \Xi_j}B_t( \cdot, \cdot ) X, 
\Xi \right\rangle \frac{\partial}{\partial X_j}G(X, \Xi)\right] = \csh^{\frac{3}{2}}\frac{|X|^2}{1 + |X|^2}\O(|\Xi|),
\end{eqnarray*}
and similarly for \eqref{hpg-err2}.  Adding these to \eqref{hpg-main-term}, we get
\begin{eqnarray}
H_{\widetilde{p}}G  & = &  \csh \left[ \frac{\langle B_t X, X \rangle}{1 + |X|^2}\right]\left( 1 + \csh^{\frac{1}{2}}
\O ( |\Xi|) \right) \label{hpg-err-2-1}\\ 
& & \quad +\csh \left[\frac{\langle B_t \Xi, \Xi \rangle}{1 + |\Xi|^2} \right]\left( 1 +\csh^{\frac{1}{2}} 
\O ( |X|)\right). \label{hpg-err-2-2}
\end{eqnarray}
Now we expand $B_t$ in a Taylor series about $(0,0)$ to get
\begin{eqnarray*}
\lefteqn{H_{\widetilde{p}}G   = } \\
& = &   \csh \left[ \frac{\langle B_t(0,0) X, X \rangle}{1 + |X|^2} + \csh^{\half} \frac{|X|^2}{1+|X|^2}
 \O(|(X, \Xi)|)\right] \cdot \\
&& \quad \quad \quad \cdot \left( 1 + \csh^{\half}\O ( |\Xi|) \right) \\ 
& & +\csh \left[\frac{\langle B_t(0,0) \Xi, \Xi \rangle}{1 + |\Xi|^2}  + \csh^{\half} \frac{|\Xi|^2}{1+|\Xi|^2} 
\O(|(X, \Xi)|)\right]\cdot \\
&& \quad \quad \quad \cdot \left( 1 +\csh^{\half} \O ( |X|)\right), 
\end{eqnarray*}
which can again be written as (\ref{hpg-err-2-1}-\ref{hpg-err-2-2}).  Recalling that $B_t(0,0)$ is block 
diagonal of the form (\ref{B-condition.a}-\ref{B-condition.b}), we get for $|(X, \Xi)| \leq \hsc^{\half} \epsilon$, 
\begin{eqnarray}
\lefteqn{ H_{\widetilde{p}} G(X, \Xi) = } \nonumber \\
& = & \left[ \sum_{j=1}^{n_{hc}} \Re \lambda_j \left( \frac{ X_{2j}^2 + X_{2j-1}^2 }{1 + |X|^2} 
+ \frac{ \Xi_{2j}^2 + \Xi_{2j-1}^2 }{1 + |\Xi|^2} \right)\right]  \left( 1 + \tilde{h}^{-\half} \O ( \epsilon )\right)  
\label{HpG-2.a} \\
& & + \left[ \sum_{j = 2n_{hc}+1}^{2n_{hc}+n_{hr}} \lambda_j \left( \frac{X_j^2}{1 + |X|^2}
 + \frac{\Xi_j^2}{1+ |\Xi|^2} \right)\right] \left( 1 + \tilde{h}^{-\half} \O ( \epsilon) \right) . \label{HpG-2.b}
\end{eqnarray}
Thus
\begin{eqnarray}
\label{Pth2}
\widetilde{P}_s(h) = \widetilde{P}(h) - ish (A( 1 + E_0))^w  + s E_1^w   + s^2 E_2^w ,
\end{eqnarray}
\noindent with $E_0 = \O(\tilde{h}^{-\half} \epsilon)$, $E_1=\O
(h^{3/2} \tilde{h}^{3/2})$, $E_2 = \O(h \tilde{h})$, and 
$A^w = \Op_{\tilde{h}}^w(A)$ for 
\begin{eqnarray}
\lefteqn{ A(X, \Xi) := } \nonumber \\ 
& = & \sum_{j=1}^{n_{hc}} \Re \lambda_j \left( \frac{ X_{2j}^2 + X_{2j-1}^2 }{1 + |X|^2} 
+ \frac{ \Xi_{2j}^2 + \Xi_{2j-1}^2 }{1 + |\Xi|^2} \right) \label{A-1-hc} \\
& & + \sum_{j=2n_{hc}+1}^{2n_{hc} + n_{hr}} \lambda_j \left( \frac{X_j^2}{1+|X|^2} 
+ \frac{\Xi_j^2}{1+ |\Xi|^2} \right). \label{A-1-hr} 
\end{eqnarray}

We claim that for $\tilde{h}$ sufficiently small, 
\begin{eqnarray}
\label{A-lower-bound}
\langle A^wU, U \rangle \geq \frac{\tilde{h}}{C} \|U \|^2
\end{eqnarray}
for some constant $C>0$, which is essentially the lower bound for the harmonic oscillator 
$\tilde{h}^2 D_X^2 + X^2$.  Clearly it suffices to prove this inequality for individual $j$ for 
the real hyperbolic terms (\ref{A-1-hr}), and in pairs for the complex
hyperbolic terms (\ref{A-1-hc}), which is the content of Lemma
\ref{harm-osc}.

Now fix $\tilde{h}>0$ and $|s|>0$ sufficiently small so that the
 estimate \eqref{A-lower-bound} holds and the errors $E_1$ and $E_2$ satisfy
\begin{eqnarray*}
\|shA^w U \|_{L^2} \gg \| s E_1^w U \|_{L^2} + \|s^2 E_2^w U \|_{L^2},
\end{eqnarray*}
and fix $\epsilon>0$ sufficiently small that the error $|E_0| \ll 1$,
independent of $h>0$. 

We now have for smooth $U$ satisfying $\Op_{\tilde{h}}^w(\psi_0(h^{\half} \bullet)) U = U + \O(h^\infty)$,
\begin{eqnarray}
\label{hsc-est}
-\Im \langle \widetilde{P}_s(h) U, U \rangle \geq \frac{h \tilde{h}}{C} \| U \|^2.
\end{eqnarray}
Now define the operator $K_h^w$ by $e^{sK_h^w} := T_{h, \tilde{h}}^{-1}
e^{sG_{\tilde{h}}^w}T_{h, \tilde{h}}$.  Translating back 
into original coordinates, and with $z \in [-1,1] + i(-c_0h + \infty)$
for sufficiently small $c_0>0$, (\ref{hsc-est}) gives 
\begin{eqnarray*}
-\Im \left\langle e^{sK_h^w} \left( P(h) - z \right) e^{-sK_h^w}u, u \right\rangle \geq \frac{h}{C_1} \| u \|^2.
\end{eqnarray*}
Finally, since $\|\exp ( \pm sK_h^w ) \| = \O (h^{-N})$ for some $N$, the theorem follows in the case where $p$ is of 
the form \eqref{assume-p}. \\
\indent For general $p$, by Proposition \ref{normal-prop-2}, there is a symplectomorphism $\kappa$ so that up to an 
elliptic factor, $\kappa^*p$ 
is of the form \eqref{assume-p}.  Using Theorem \ref{gamma-egorov} to quantize $\kappa$ as an $h$-FIO
$F$, we get
\begin{eqnarray*}
\Op_h^w \left( \kappa^* p + E_1 \right) = F^{-1} P(h) F ,
\end{eqnarray*}
where $E_1 = \O(h^2)$ is the error arising from Theorem
\ref{gamma-egorov}.  We may then 
use the previous argument for $\kappa^*p$ getting an additional error of $\O(h^2)$ from Theorem \ref{gamma-egorov} in 
\eqref{hsc-est}, which is the same order as $E_1$.
\end{proof}

\begin{remark}
The error arising at the end of the proof of Theorem
\ref{main-theorem-1} from the use of Theorem \ref{gamma-egorov} is of
order $\O(h^2)$ and hence negligible compared to our lower bound of
$h$ for $A$.  However, the estimate of $A$ is used for the imaginary
part of $\widetilde{P}_s$, and the error in Theorem \ref{gamma-egorov}
is real, so $\O(h)$ would
have been sufficient.
\end{remark}

\begin{lemma}
\label{harm-osc}
Let
\begin{eqnarray*}
a_0(y, \eta):= \frac{y_j^2}{ \langle y \rangle^2} + \frac{\eta_j^2}{ \langle \eta \rangle^2},
\end{eqnarray*}
for $(y, \eta) \in \reals^{2n-2}$, and $\langle y \rangle = ( 1 +
|y|^2)^{1/2}$, and let 
\begin{eqnarray*}
a_1(y, \eta) := \frac{y_{2j}^2 + y_{2j-1}^2}{\langle y \rangle^2} +
\frac{\eta_{2j}^2 + \eta_{2j-1}^2}{\langle \eta \rangle^2}.
\end{eqnarray*}
Then $a_i$, $i=0,1$ satisfies
\begin{eqnarray}
\label{ai-ineq}
\langle \Op_{\tilde{h}}^w(a_i)U, U \rangle \geq \frac{\tilde{h}}{C} \|U \|^2
\end{eqnarray}
for $\tilde{h}>0$ sufficiently small and a constant $0<C<\infty$.
\end{lemma}
\begin{proof}
The idea of the proof is that $a_i$ is essentially the harmonic
oscillator which satisfies the inequality \eqref{ai-ineq}.  We write each $a_i$ as a $a_i = |b|^2$ for $b$ a complex
symbol.  Observe $a_0(y, \eta) = |b(y, \eta)|^2$ with 
\begin{eqnarray*}
b(y, \eta): = \frac{y_j}{\langle y \rangle} + i \frac{\eta_j}{\langle \eta \rangle}.
\end{eqnarray*}
Thus, using the $\tilde{h}$-Weyl calculus,
\begin{eqnarray}
\label{A_0}
a_0^w(y, \tilde{h}D_y) = b^w(y, \tilde{h}D_y)^* b^w(y, \tilde{h} D_y) + c^w(y, \tilde{h} D_y),
\end{eqnarray}
where 
\begin{eqnarray}
c(y, \eta) & = & \tilde{h} \left\{  \frac{\eta_j}{\langle \eta \rangle},
\frac{y_j}{\langle y \rangle}\right\} 
+ \O ( \tilde{h}^2) \nonumber \\
 & = & \tilde{h} \langle y \rangle^{-3} \langle \eta \rangle^{-3}\left( 1 + \O(|y|^2 + |\eta|^2) \right) 
+ \O( |y|^2|\eta|^2) + \O ( \tilde{h}^2)  \label{C}. 
\end{eqnarray}
For $(y, \eta)$ small, $c$ is bounded from below by $\tilde{h}$ as in
(\ref{A-lower-bound}), and for large $(y, \eta)$ we have 
\be
C^{-1} \leq a_0 \leq C
\ee
for some constant $C>0$.  Hence for large $(y, \eta)$, \eqref{A_0} is
bounded from below independent of $\tilde{h}$.  Observe $a_1(y, \eta)
= |b_{2j}(y, \eta)|^2 + |b_{2j-1}(y, \eta)|^2$ for 
\begin{eqnarray*}
b_k(y, \eta) = \frac{y_k}{\langle y \rangle } -i \frac{\eta_k}{\langle \eta \rangle },
\end{eqnarray*}
and the same argument applies to $a_1$ as to $a_0$. 
\end{proof}

\begin{remark}
It is interesting to note that the estimate (\ref{main-theorem-1-est}) depends only on the real parts of the eigenvalues 
$\lambda_j$ above.  Unraveling the definitions, the eigenvalues $\lambda_j$ are logarithms of the eigenvalues of the 
linearized Poincar\'{e} map $dS(0)$ from above.  Then (\ref{main-theorem-1-est}) depends only on the modulis of the 
eigenvalues of $dS(0)$.
\end{remark}

\numberwithin{equation}{subsection}
\subsection{A Return to Quadratic Forms}
Recall the only place we have used that the eigenvalues are distinct
is in determining the possible form of the quadratic form $q(\rho)$
defined by $dS(0) = \exp H_q$.  We then considered the Hamilton, or
Fundamental matrix $B$ defined by
\begin{eqnarray}
\label{qB}
q(\rho) = : \half \tom (\rho, B \rho).
\end{eqnarray}
We follow \cite{hormander5} and return to the setup for Proposition
\ref{normal-prop-1}.  All of the following changes of variables will
be linear, so we may assume we are working in $\reals^{2n-2}$ and
choose local symplectic coordinates in
which $\tom$ is the standard symplectic form
\begin{eqnarray*}
\tom = \sum_{j=1}^{n-1} d\xi_j \wedge dx_j.
\end{eqnarray*}
Then we can write \eqref{qB} in a more easily manipulated form:
\begin{eqnarray*}
q(\rho) = : \half \langle \rho , JB \rho \rangle
\end{eqnarray*}
where $J$ is the matrix of symplectic structure on $\reals^{2n}$,
\begin{eqnarray*}
J = \left( \begin{array}{cc} 0 & -I \\ I & 0 \end{array} \right).
\end{eqnarray*}
As mentioned previously, the eigenvalues of $B$ are the logarithms of
the eigenvalues of $dS(0)$ (with a suitably chosen branch cut), hence
have nonzero real part, and come in pairs $\lambda, - \lambda$ for the
positive real
hyperbolic eigenvalues, and $4$-tuples $\lambda, - \lambda,
\bar{\lambda}, - \bar{\lambda}$ for the complex hyperbolic.  If we allow $\rho$ to be complex for the moment, and
denote by $V_\lambda$ the generalized eigenspace for $\lambda$ real or
complex, we see 
\begin{eqnarray*}
\tom (V_\lambda, V_{\lambda'}) = 0
\end{eqnarray*}
unless $\lambda + \lambda' = 0$.  We then consider the spaces
$V_\lambda \oplus V_{-\lambda}$, which is symplectic with the
restricted symplectic form $\left. \tom \right|_{V_\lambda \oplus
  V_{-\lambda}}$, since $\lambda \neq 0$.  As in \S
\ref{symplectic-geometry} we choose the pairs and $4$-tuples of
eigenvalues so that $\Re \lambda >0$ and $\Im \lambda \geq 0$.  We thus have a
decomposition of $\reals^{2n-2}$ into symplectic subspaces
\begin{eqnarray*}
\reals^{2n-2} =  \left( \bigoplus_{j = 1}^{n_{hc}} V_{\lambda_j}
\oplus V_{-\lambda_j} \oplus V_{\bar{\lambda}_j}
\oplus V_{-\bar{\lambda}_j} \right) \bigoplus \left( \bigoplus_{j =
  n_{hc}+1}^{n_{hc} + n_{hr}}  V_{\lambda_j}
\oplus V_{-\lambda_j}\right)
\end{eqnarray*}
where $n_{hr}$ is the number of real eigenvalues with $\lambda>0$ and $n_{hc}$ is the
number of complex eigenvalues with $\Re \lambda >0, \Im \lambda >0$.  Our notation here means if $\lambda_j$
has multiplicity $k_j$, then
\begin{eqnarray*}
\sum_{j=1}^{n_{hc}} 4k_j + \sum_{j=n_{hc} + 1}^{n_{hc} + n_{hr}} 4 k_j = 2n-2.
\end{eqnarray*}
Fix $\lambda$ real or complex, $\Re \lambda >0$, $\Im \lambda \geq 0$,
with multiplicity
greater than $1$ and consider the
complex symplectic subspace $V_\lambda
\oplus V_{-\lambda}$.  Assume $V_\lambda$ has dimension $m$.  Note $B$
restricts to a linear map in $V_\lambda$,
$T:= \left. B \right|_{V_{\lambda}}$, such that $T - \lambda I$ is
nilpotent.  Our definitions equip $V_\lambda \oplus V_{-\lambda}$
with a symplectic structure in which $V_{-\lambda}$ is dual and
isomorphic to $V_{\lambda}$.  We abuse notation and write a point
$(x, \xi) \in V_\lambda \oplus V_{-\lambda}$.  Then if we put $T$ into Jordan form in
$V_\lambda$ so that $Tx = \lambda x + (x_2, x_3, \ldots, x_m, 0)$, we
obtain a symplectic change of coordinates by writing
\begin{eqnarray*}
\left. B \right|_{V_\lambda \oplus V_{-\lambda}}(x, \xi) = (\lambda x
+ (x_2, \ldots, x_m, 0), -\lambda \xi - (0, \xi_1, \xi_2, \ldots,
\xi_{m-1})),
\end{eqnarray*}
by the symplectic skew symmetry of $B$.  In these coordinates we then have
$q_\lambda$, the
projection of $q$ onto $V_\lambda \oplus V_{-\lambda}$,
\begin{eqnarray}
\label{q-form-cx}
q_\lambda (x,\xi) = \lambda \sum_{l=1}^k x_l \xi_l + \sum_{l=1}^{k-1}
x_{l+1}\xi_l,
\end{eqnarray}
where $k$ is the multiplicity of $\lambda$.  This is the normal form
in complex variables, with the ``actions'' $\lambda x_j \xi_j$ as in
\S \ref{symplectic-geometry}, but with the additional terms coming from the Jordan
form.  In order to understand the real normal form, there are two cases to
examine.

{\bf Case 1: $\lambda>0$ is real.}  Then the space $V_\lambda \oplus
V_{-\lambda}$ is real, the change of variables above is real, and we
get $q_\lambda$ exactly as in \eqref{q-form-cx}.  Let the real matrix
$Q_\lambda$ be defined by the real normal form:
\begin{eqnarray}
\label{Q-matrix}
q_\lambda(x,\xi) = : \half \langle (x, \xi), Q (x, \xi) \rangle.
\end{eqnarray}
Then $Q$ takes the special form
\begin{eqnarray*}
Q = \left( \begin{array}{cc} 0 & A  \\ A^T & 0 \end{array} \right) 
\end{eqnarray*}
where $A$ is the $k \times k$ matrix
\begin{eqnarray}
\label{A-1}
A = \left( \begin{array}{cccc} \lambda & 0 &
  \multicolumn{2}{c}\dotfill \\ 1 & \lambda & 0 & \ldots \\
0 & \ddots & \ddots & 0 \\
\vdots & \ldots & 1 & \lambda
\end{array}
\right)
\end{eqnarray}
and $A^T$ denotes the transpose of $A$.

{\bf Case 2: $\lambda$ complex, $\Re \lambda >0$, $\Im \lambda >0$.}  We use a similar 
change of variables to that in \S \ref{symplectic-geometry}.  That is, let
$\{e_l, f_l\}$ be the generalized eigenvectors for $\lambda, -\lambda$
respectively.  Here, $1 \leq l \leq k$ where $k$ is the multiplicity
of $\lambda$.  Then
$\{e_l,f_l,\bar{e}_l, \bar{f}_l \}$ forms a basis for a complex vector space which is
the complexification of a real symplectic vector space.  We then
consider the projection $q_\lambda$ of $q$ onto the space 
\begin{eqnarray*}
W = V_\lambda
\oplus V_{-\lambda} \oplus V_{\bar{\lambda}} \oplus
V_{-\bar{\lambda}}.
\end{eqnarray*}
Write a point $\rho$ in $W$ as 
\begin{eqnarray*}
\rho = \sum_{l=1}^k z_l e_l + \zeta_l f_l + w_l \bar{e}_l + \eta_l \bar{f}_l,
\end{eqnarray*}
so that 
\begin{eqnarray*}
q_\lambda(\rho) = \lambda \sum_1^k z_l \zeta_l + \bar{\lambda}
\sum_1^k w_l \eta_l + \sum_1^{k-1} z_{l+1} \zeta_l + \sum_1^{k-1}
w_{l+1} \eta_l.
\end{eqnarray*}
We define as in \S \ref{symplectic-geometry} a real symplectic basis
$\{e_l^1, e_l^2, f_l^1, f_l^2\}$ for $1 \leq l \leq k$ by 
\begin{eqnarray*}
e_l = \frac{1}{\sqrt{2}}(e_l^1 + ie_l^2), \quad f_l= \frac{1}{\sqrt{2}} (f_l^1 -
i f_l^2),
\end{eqnarray*}
and write in these new coordinates
\begin{eqnarray*}
\rho = \sum_{l=1}^k \sum_{r=1}^2 x_l^r e_l^r + \xi_l^r r_l^r.
\end{eqnarray*}
Then we get the real normal form of $q_\lambda$ in these coordinates:
\begin{eqnarray*}
q_\lambda(\rho) & = & \Re \lambda \sum_1^k \left( x_{2l-1} \xi_{2l-1} +
x_{2l} \xi_{2l} \right) - \Im \lambda \sum_1^k \left( x_{2l}\xi_{2l-1}
- x_{2l-1} \xi_{2l} \right)  \\
&& \quad \quad + \sum_1^{k-1} \left( x_{2l+1}\xi_{2l-1} + x_{2l+2}
\xi_{2l} \right).
\end{eqnarray*}
We again define the real matrix $Q$ in terms of the real quadratic
normal form $q_\lambda$ by \eqref{Q-matrix}, which now takes the form
\begin{eqnarray*}
Q = \left( \begin{array}{cc} 0 & A \\ A^T & 0 \end{array} \right),
\end{eqnarray*}
where $A$ is the $2k \times 2k$ matrix
\begin{eqnarray} 
\label{A-2}
\left( \begin{array}{cccc} 
\Lambda & 0 &
  \multicolumn{2}{c}\dotfill \\ I & \Lambda & 0 & \ldots \\
0 & \ddots & \ddots & 0 \\
\vdots & \ldots & I & \Lambda
\end{array} \right),
\end{eqnarray}
with $I$ the $2 \times 2$ identity matrix and 
\begin{eqnarray*}
\Lambda = \left( \begin{array}{cc}
\Re \lambda & - \Im \lambda \\
\Im \lambda & \Re \lambda
\end{array} \right).
\end{eqnarray*}
Putting this discussion together with the proof of Proposition
\ref{normal-prop-1}, we have proved the following:
\begin{proposition}
\label{normal-prop-3}
Let $p \in \Ci ( T^*X)$, $\gamma \subset \{p = 0 \}$ as above, with
the linearized Poincar\'{e} map $dS(0)$ having eigenvalues $\{\mu_j\}$ not on the unit
circle, and suppose $\mu_j$ has multiplicity $k_j$.  Then there exists a neighbourhood, $U$, of $\gamma$ in $T^*X$, 
a smooth positive function $b \geq C^{-1} >0$ defined in $U$, and a 
symplectomorphism $\kappa : U \to \kappa(U) \subset T^*\SS_{(t, \tau)}^1 \times T^* \reals_{(x, \xi)}^{n-1}$ such that 
\begin{eqnarray*}
\kappa(\gamma)   =  \{ (t,0;0,0) : t \in \SS^1 \}, 
\end{eqnarray*}
and $b(t, \tau, x, \xi) p = \kappa^*(g + r)$, with
\begin{eqnarray*}
\lefteqn{ g(t, \tau; x, \xi) =} \nonumber \\ 
& = &  \tau + \sum_{j=1}^{n_{hc}} \sum_{l=1}^{k_j}\left( \Re \lambda_j \left( x_{2l-1} \xi_{2l-1} + x_{2l} \xi_{2l} \right) 
- \Im \lambda_j \left( x_{2l-1} \xi_{2l} - x_{2l}\xi_{2l-1} \right)
\right)  \\
&& \quad \quad + \sum_{j=1}^{n_{hc}} \sum_{l=1}^{k_j-1} \left( x_{2l+1}\xi_{2l-1} + x_{2l+2}
\xi_{2l} \right) \\
 & & \quad \quad +  \sum_{j = 2n_{hc}+1}^{2n_{hc} + n_{hr}} \left(\sum_{l=1}^{k_j}
\lambda_j x_l \xi_l  + 
\sum_{l=1}^{k_j - 1} x_{l+1}\xi_l\right),
\end{eqnarray*}
where $\lambda_j = \log \mu_j$ for each $j$ (with a suitable branch cut) and $r  =  \O (|x|^3 + |\xi|^3)$. 
\end{proposition}

The proof of Lemma \ref{zero-cov} depends only on the moduli of the
eigenvalues of $dS(0)$ restricted to the stable and unstable
manifolds, hence does not depend on the multiplicities, or the Jordan
form.  Consequently we have the analogue of Proposition \ref{normal-prop-2}.

\begin{proposition}
\label{normal-prop-4}
Under the assumptions of Proposition \ref{normal-prop-3}, there exists
a neighbourhood, $U$, of $\gamma$ in $T^*X$, a smooth
positive function $b\geq C^{-1} >0$ defined in $U$, 
a symplectomorphism 
$\kappa : U \to \kappa(U) \subset T^*\SS_{(t, \tau)}^1 \times T^* \reals_{(x, \xi)}^{n-1}$, 
and a smooth, $n \times n$-matrix valued function $B_t$ such that 
\begin{eqnarray*}
\kappa(\gamma)  & = & \{ (t,0;0,0) : t \in \SS^1 \}, \quad
\mathrm{and} \,\,\, b(t, \tau, x, \xi) p = \kappa^*g, \,\,\, 
\mathrm{with} \\
g(t, \tau; x, \xi) & = &\tau + \langle B_t(x, \xi) x, \xi \rangle,
\end{eqnarray*}
with $B_t$ satisfying 
\begin{eqnarray*}
\lefteqn{\left\langle B_t(0,0)x, \xi \right\rangle = }  \\ 
& = &  \sum_{j=1}^{n_{hc}} \sum_{l=1}^{k_j}\left( \Re \lambda_j \left( x_{2l-1} \xi_{2l-1} + x_{2l} \xi_{2l} \right) 
- \Im \lambda_j \left( x_{2l-1} \xi_{2l} - x_{2l}\xi_{2l-1} \right)
\right)  \\
&& \quad \quad + \sum_{j=1}^{n_{hc}} \sum_{l=1}^{k_j-1} \left( x_{2l+1}\xi_{2l-1} + x_{2l+2}
\xi_{2l} \right) \\
 & & \quad \quad +  \sum_{j = 2n_{hc}+1}^{2n_{hc} + n_{hr}} \left( \sum_{l=1}^{k_j}
\lambda_j x_l \xi_l  + \sum_{l=1}^{k_j - 1} x_{l+1}\xi_l \right).
\end{eqnarray*}
\end{proposition}

\subsection{End of the Proof of Theorem \ref{main-theorem-1}}
\label{non-distinct-section}
Now we turn our attention to the proof of Theorem
\ref{main-theorem-1} in the case of non-distinct eigenvalues
of $dS(0)$.  Recall the key feature to the proof of Theorem
\ref{main-theorem-1} in the case of distinct eigenvalues was that the normal form given in Proposition
\ref{normal-prop-2} has quadratic part $q(x, \xi)$ with the
property that there exists another quadratic form 
\begin{eqnarray*}
w(x, \xi) = \half \langle W(x, \xi), (x, \xi) \rangle
\end{eqnarray*}
such that $H_q w(x, \xi)$ is a positive definite quadratic form.  Then
we would like our escape function to be $G(x, \xi) = w(x, \xi)$,
however for technical reasons we
had to use a logarithmic escape function and form the families $e^{\pm
  s G^w}$.  With the following theorem, the proof of Theorem
\ref{main-theorem-1} is complete.
\begin{theorem}
\label{non-distinct-theorem}
Suppose $q \in \Ci (\reals^{2m})$ is quadratic of the form
\begin{eqnarray}
\lefteqn{ q(x,\xi)  = } \nonumber \\ 
&&  =\sum_{j=1}^{n_{hc}} \sum_{l=1}^{k_j}\left( \Re \lambda_j \left( x_{2l-1} \xi_{2l-1} + x_{2l} \xi_{2l} \right) 
- \Im \lambda_j \left( x_{2l-1} \xi_{2l} - x_{2l}\xi_{2l-1} \right)
\right)  \label{p-1a}\\
&& \quad \quad + \sum_{j=1}^{n_{hc}} \sum_{l=1}^{k_j-1} \left( x_{2l+1}\xi_{2l-1} + x_{2l+2}
\xi_{2l} \right) \label{p-1b}\\
 & & \quad \quad +  \sum_{j = 2n_{hc}+1}^{2n_{hc} + n_{hr}} \left( \sum_{l=1}^{k_j}
\lambda_j x_l \xi_l  + \sum_{l=1}^{k_j - 1} x_{l+1}\xi_l \right), \label{p-1c}
\end{eqnarray}
and 
\begin{eqnarray*}
G(x, \xi) = \half \left( \log (1 + |x|^2) - \log(1 + |\xi|^2 ) \right).
\end{eqnarray*}
Then there exist $m \times m$ nonsingular matrices $A$
and $A'$, positive real numbers $0 < r_1 \leq r_2, \leq \cdots \leq
r_{m}< \infty$, and symplectic coordinates $(x, \xi)$ such that
\begin{eqnarray}
\label{non-distinct-theorem-statement}
H_q(G) = \frac{\sum_{j=1}^m r_j^{-2} x_j^2 }{1 + |A x|^2} +
\frac{\sum_{j=1}^m r_j^{-2} \xi_j^2 }{1 + |A' \xi|^2}.
\end{eqnarray} 
\end{theorem}

\begin{proof}
First, suppose
\begin{eqnarray*}
g(x, \xi) = \half \left\langle \tilde{g} \left( \begin{array}{c} x \\ \xi
\end{array} \right), \left( \begin{array}{c} x \\ \xi
\end{array} \right) \right\rangle
\end{eqnarray*}
is a real quadratic form with $\tilde{g}$ symmetric of the form
\begin{eqnarray*}
\tilde{g} = \left( \begin{array}{cc} P & 0 \\ 0 & -P \end{array} \right),
\end{eqnarray*}
where $P$ is symmetric and nonsingular.  Then 
\begin{eqnarray*}
\partial_{x}  \half \log( 1 + \langle Px,x \rangle )  =  \frac{Px}{1 +
  \langle Px,x \rangle},
\end{eqnarray*}
and similarly for $\xi$ so studying
\begin{eqnarray*}
H_q \left( \half \left( \log(1 + \langle Px,x \rangle) - \log(1 + \langle
P\xi,\xi \rangle ) \right) \right)
\end{eqnarray*}
can be reduced to studying $H_q g(x, \xi)$, modulo the positive terms
$1 + \langle P \cdot , \cdot \rangle$ in the denominator.  If $q(x, \xi)$ is of the
form (\ref{p-1a}-\ref{p-1c}), then we can write $q$ in terms of the
fundamental matrix $B$:
\begin{eqnarray*}
q(x, \xi) = \left\langle 
\left( \begin{array}{c} x \\ \xi
\end{array} \right), JB \left( \begin{array}{c} x \\ \xi
\end{array} \right) \right\rangle,
\end{eqnarray*}
where 
\begin{eqnarray*}
J = \left( \begin{array}{cc} 0 & -I \\ I & 0 \end{array} \right)
\end{eqnarray*}
as usual.  Then the vector field $H_q$ can be written as
\begin{eqnarray*}
H_q = \left\langle B \left( \begin{array}{c} x \\ \xi
\end{array} \right), \left( \begin{array}{c} \partial_x \\ \partial_\xi
\end{array} \right) \right\rangle,
\end{eqnarray*}
and
\begin{eqnarray*}
H_qg & = & \left\langle B \left( \begin{array}{c} x \\ \xi
\end{array} \right), \left( \begin{array}{c} \partial_x \\ \partial_\xi
\end{array} \right) \right\rangle \left( \half \left\langle \tilde{g} \left( \begin{array}{c} x \\ \xi
\end{array} \right), \left( \begin{array}{c} x \\ \xi
\end{array} \right) \right\rangle \right) \\
& = & \left\langle B \left( \begin{array}{c} x \\ \xi
\end{array} \right), \tilde{g} \left( \begin{array}{c} x \\ \xi
\end{array} \right) \right\rangle \\
& = & \left\langle \tilde{g}B \left( \begin{array}{c} x \\ \xi
\end{array} \right),  \left( \begin{array}{c} x \\ \xi
\end{array} \right) \right\rangle,
\end{eqnarray*}
since $\tilde{g}$ is symmetric.

Now from the discussion preceding the statement of Theorem
\ref{non-distinct-theorem}, we know $B = -JQ$ for $Q$ of the form
\begin{eqnarray}
\label{Q-full}
Q = \left( \begin{array}{cc} 0 & A \\ A^T & 0 \end{array} \right),
\end{eqnarray}
where $A$ is block diagonal with diagonal elements of the form
\eqref{A-1} or \eqref{A-2}.  Thus with the same $A$ as \eqref{Q-full},
\begin{eqnarray*}
B = \left( \begin{array}{cc} A^T & 0 \\ 0 & -A \end{array} \right).
\end{eqnarray*}
Now we have reduced the problem to finding nonsingular $P$ such that $PA^T$ and
$PA$ are both positive definite.  But we know that if $\lambda$ is an
eigenvalue of $A$, then $\Re \lambda >0$, so $A$ is positive definite
and $P = I$ suffices.  \eqref{non-distinct-theorem-statement} then
follows immediately from Lemma \ref{HZ-ell-lemma}.
\end{proof}
We have used the following classical lemma (see, for example,
\cite{hz} for a proof).

\begin{lemma}
\label{HZ-ell-lemma}
Let 
\begin{eqnarray*}
q(x, \xi ) = \half \langle Q (x, \xi), (x, \xi) \rangle
\end{eqnarray*}
be a positive definite quadratic form, where $Q$ is symmetric.  Then
there are positive numbers $0 < r_1 \leq r_2 \leq \cdots \leq r_{m}
< \infty$ and a linear symplectic transformation $T$ such that 
\begin{eqnarray*}
q(T(x, \xi)) = \sum_{j=1}^{m} \frac{1}{r_j^2} (x_j^2 + \xi_j^2 ).
\end{eqnarray*}
Further, if $T'$ is another linear symplectic transformation such that
\begin{eqnarray*}
q(T'(x, \xi)) = \sum_{j=1}^{m} \frac{1}{r_j'^2} (x_j^2 + \xi_j^2 )
\end{eqnarray*}
for $0 < r_1' \leq \cdots \leq r_{m}' < \infty$, then $r_j = r_j'$ for
all $j$ and $T = T'$.
\end{lemma}


\section{Proof of Theorem \ref{main-theorem-2} and the Main Theorem}
\label{main-theorem-2-proof} 
\subsection{Proof of Theorem \ref{main-theorem-2}}
\begin{figure}
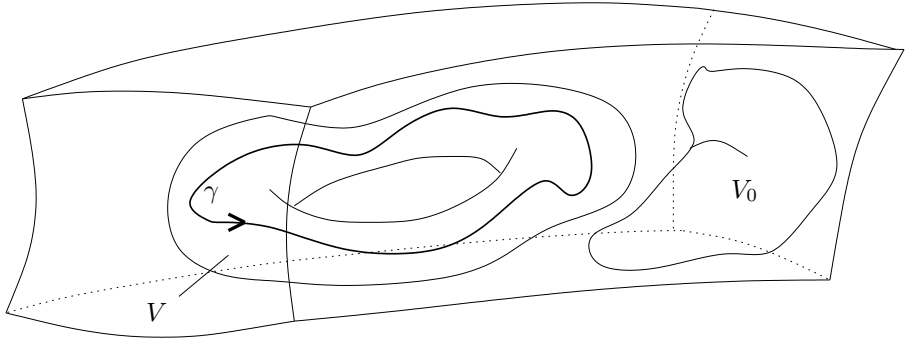

\include{fig2}
\caption{\label{fig-p0} The energy surface $\{ p^{-1}(0) \}$.}
\end{figure}

In this section we show how to use Theorem \ref{main-theorem-1} with a few other results to deduce 
Theorem \ref{main-theorem-2}.  This is similar to \cite{BZ}, with the generalization of the 
loxodromic assumption.  First we need the following standard lemma.  
\begin{lemma}
\label{wf-lemma}
Suppose $V_0 \Subset T^*X$, $p$ is a symbol, $T >0$, $A$ an operator, and $V \Subset T^*X$ a neighbourhood of $\gamma$ satisfying
\begin{eqnarray}
\left\{ \begin{array}{l}
\forall \rho \in \{ p^{-1}(0) \} \setminus V, \,\,\, \exists \, 0 <t<T \,\, \text{and} \,\, \epsilon = \pm 1 \,\,\, 
\text{such that} \\
\exp(\epsilon s H_p)(\rho) \subset \{ p^{-1}(0) \} \setminus V \,\, \text{for} \,\, 0 < s < t, \,\, \text{and} \\
\exp(\epsilon t H_p)(\rho) \in V_0; \\
\end{array} \right.
\end{eqnarray}
and $A$ is microlocally elliptic in $V_0 \times V_0 $.  If $B \in \Psi^{0,0}(X, \Omega_X^\half)$ and $\WF (B) \subset T^*X 
\setminus V$, then
\begin{eqnarray*}
\left\| Bu \right\| \leq C \left( h^{-1} \left\| Pu \right\| + \| Au \| \right) + \O (h^\infty) \|u\|.
\end{eqnarray*}
\end{lemma}
Figure \ref{fig-p0} is a picture of the setup of Lemma \ref{wf-lemma}.
\begin{proof}
Since $\{ p^{-1}(0) \}$ is compact, we can replace $V_0$ with a precompact neighbourhood of $V_0 \cap \{p^{-1}(0) \}$.  
We will prove a local version which can be pasted together to get the global estimate.  We may assume 
$\WF(A) \subset U$, where $U$ is a small open neighbourhood of some point $\rho_0 \in V_0$, and 
\begin{eqnarray}
\WF (B) \Subset \bigcup_{0 \leq t \leq t_0} \exp (\epsilon t H_p) ( U_1 ) \subset T^*X \setminus V,
\end{eqnarray}
where $U_1 \Subset U$ and $A$ is microlocally elliptic on $U_1\times U_1$.  For $|t|\leq t_1$ 
sufficiently small, by Proposition \ref{hDx-prop} there is a 
microlocally invertible $h$-FIO $T$ which   
conjugates $P$ to $hD_{x_1}$.  Set $\tilde{u} = Tu$, and let
$\tilde{B} \in \Psi^{0,0}$ be microlocally $1$ on $\WF (B) \times \WF (B)$ 
and $0$ microlocally outside $\left( \cup_{0 \leq t \leq t_1} \exp
(\epsilon t H_p) ( U_1 ) \right)^2 \subset (T^*X \setminus V)^2$.  
We calculate 
\begin{eqnarray*}
\frac{1}{2} \partial_{x_1} \|\tilde{u} \|^2 & = & \langle \partial_{x_1} \tilde{u}, \tilde{u} \rangle \\
& \leq & \left\| \partial_{x_1} \tilde{u} \right\| \|\tilde{u}\| \\
& \leq & \frac{1}{4} h^{-1} \left\| TPT^{-1}\tilde{u} \right\|^2 + \|\tilde{u} \|^2 + \O(h^\infty) \|u\|_{L^2(X)}^2 \\
\implies \| \tilde{B}T^{-1}\tilde{u} \|^2_{L^2(X)} & \leq & C_{t_1} \Big( h^{-1} \left\| TPT^{-1}\tilde{u} \right\|^2_{L^2(X)} 
+ \left\|A T^{-1}\tilde{u} \right\|^2_{L^2(X)} + \\
&& \quad \quad \quad + \O(h^\infty)\|u\|^2_{L^2(X)}\Big),
\end{eqnarray*}
where the last inequality follows from Gronwall's inequality.  But $\|BT^{-1}\tilde{u}\|^2_{L^2(X)} \leq 
\|\tilde{B}T^{-1}\tilde{u} \|^2_{L^2(X)}$ 
gives the result for small $t$.  Then we partition $[0, t_0]$ into finitely many subintervals and apply the small 
$t$ argument to each one.
\end{proof}
Using this lemma, we can deduce the following proposition.
\begin{proposition}
\label{Q(z)u-prop}
Suppose $\psi_0 \in \s^{0,0}(T^*X) \cap \Ci_c(T^*X)$ is a microlocal cutoff function to a small neighbourhood 
of $\gamma \subset \{p^{-1}(0) \}$.  For $Q(z) = P(h) - z - iCha^w$ as above with $z \in [-1,1] + i(-c_0h, \infty)$, 
$c_0 >0$ and $C>0$ sufficiently large, we have
\begin{eqnarray}
\label{Q(z)u-prop-est}
Q(z)u = f \Longrightarrow \left\| (1 - \psi_0 )^w u \right\| \leq C h^{-1} \|f \| + \O(h^\infty) \|u \|.
\end{eqnarray}
\end{proposition}
For this proposition and the proof, we use the convenient shorthand notation: for a symbol $b$, $b^w := \Op_h^w(b)$.
\begin{remark}
Note that Proposition \ref{Q(z)u-prop} is the best possible situation.  It says roughly that away from $\gamma$, 
$Q^{-1}$ is bounded by $Ch^{-1}$.  Thus the global statement in Theorem \ref{main-theorem-2} represents a loss of 
$\sqrt{\log (1/h)}$.
\end{remark}
\begin{proof}
Choose $c_0 >0$ from Theorem \ref{main-theorem-1}, microlocal cutoff
functions $\psi_1$, $\psi_2$ such that $\WF (1-\psi_j) \cap \gamma =
\emptyset$, and $C>0$ sufficiently large so that  
\begin{eqnarray*}
(Ca-c_0)^w (1 - \psi_1)^w \geq \left\{ \begin{array}{l}
c_0 (1 - \psi_1)^w /2\\
c_0 ((1 - \psi_2)^w)^*(1-\psi_2)^w/2, \end{array} \right.
\end{eqnarray*}
and $\supp \psi_1 \subset \{ \psi_2 = 1 \}$ (see Figure \ref{fig-cutoffs-1}).  
\begin{figure}
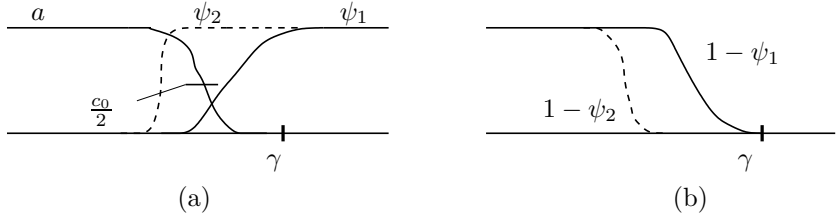

\include{fig4}
\caption{\label{fig-cutoffs-1} (a) The cutoff functions $a$, $\psi_1$, and $\psi_2$.  (b) $(1-\psi_2)^2 \leq (1-\psi_1)$.}
\end{figure}
Then we calculate
\begin{eqnarray*}
\frac{1}{2} c_0 h \int_X \left| (1 - \psi_2)^w u \right|^2 dx & \leq & h \int_X \left( Ca^w + h^{-1} \Im z \right) 
u \overline{ (1- \psi_1)^w u } dx \\
& = & -\Im \int_X Q(z) u \overline{ (1-\psi_1)^w u} dx \\
& = & -\Im \int_X f \overline{ (1-\psi_1)^w u} dx \\
& \leq & \|f\| \left( \left\| (1 - \psi_1)^w u \right\| + \O(h^\infty) \|u \| \right) \\
& \leq & (4\epsilon h)^{-1} \|f \|^2 + \epsilon h \left\| (1 -\psi_1)^w u \right\|^2 + \O(h^\infty)\|u\|^2
\end{eqnarray*}
Now we use Lemma \ref{wf-lemma} with $A = (1-\psi_2)^w$, $B = (1-\psi_1)^w$, and $P = Q(z)$, which we may do since 
the perturbation terms in $Q(z)$ are all of lower order.  Thus
\begin{eqnarray*}
\left\|(1-\psi_1)^w u \right\| & \leq & Ch^{-1} \left\|Q(z) u \right\|  + \left\| (1 - \psi_2)^w u \right\| 
+ \O(h^\infty) \|u \| \\
\implies \left\|(1-\psi_1)^w u \right\|^2 & \leq & Ch^{-1} \|f\| \left( Ch^{-1} \|f\| + \left\| (1 - \psi_2)^w 
u \right\| \right) + \\ 
&& \quad \quad + \left\|( 1 -\psi_2)^w u \right\|^2 + \O(h^\infty) \|u \|^2 \\
& \leq & Ch^{-2} \|f\|^2 + \left\|(1 - \psi_2)^w u \right\|^2 + \O(h^\infty) \|u \|^2 \\
& \leq & Ch^{-2} \|f\|^2 + \epsilon  \left\| (1 -\psi_1)^w u \right\|^2 + \O(h^\infty) \|u \|^2,
\end{eqnarray*}
which gives \eqref{Q(z)u-prop-est} with $\psi_0$ replaced by $\psi_1$.  Another application of Lemma \ref{wf-lemma} 
with $A = (1 - \psi_2)^w$, $B = (\psi_1 - \psi_0)^w$, and $P = Q(z)$ shows the error $\|(\psi_1 - \psi_0)^w u \|$ is
 bounded by the same estimate as in \eqref{Q(z)u-prop-est}.
\end{proof}
We will need the next lemma, which is essentially an operator version of the classical Three-Line Theorem from complex 
analysis.  We include the proof here for the reader's convenience,
collected from \cite{BZ}, \cite{Bur}), and \cite{TaZw}.
\begin{lemma}
\label{BZ-lemma}
Let $\mathcal{H}$ be a Hilbert space, and assume $A,B: \mathcal{H} \to \mathcal{H}$ are bounded, self-adjoint operators 
satisfying $A^2 = A$ and $BA = AB = A$.  Suppose $F(z)$ is a family of bounded operators satisfying 
$F(z)^* = F(\bar{z})$, $\Re F \geq C^{-1} \Im z$ for $\Im z >0$, and further assume
\begin{eqnarray*}
BF^{-1}(z) B \,\, \text{is holomorphic in}\,\, \Omega:= [-\epsilon, \epsilon] + i[-\delta, \delta], \,\, 
\text{for} \,\, \frac{\delta}{\epsilon} \ll M^{-\frac{1}{N_1}} < 1
\end{eqnarray*}
for some $N_1>0$, where $\|B F^{-1}(z)B \| \leq M$.  Then for $|z| < \epsilon/2$, 
$\Im z = 0$,
\begin{eqnarray*}
(a) \quad \left\| B F^{-1}(z) B \right\| & \leq & C \frac{\log M}{\delta}, \\
(b) \quad \left\| B F^{-1}(z) A \right\| & \leq & C \sqrt{ \frac{\log M}{\delta}}.
\end{eqnarray*}
\end{lemma}
\begin{proof}
For the proof of part (a), consider the holomophic operator-valued function $f(z) = B F(z)^{-1} B$.  Choose 
$\psi \in \Ci_c([-3 \epsilon/4, 3 \epsilon/4])$, $\psi \equiv 1$ on $[-\epsilon/2, \epsilon/2]$, and for 
$z \in \Omega$, set
\begin{eqnarray*}
\phi(z) = \delta^{-\half} \int e^{-(x-z)^2/\delta} \psi(x) dx.
\end{eqnarray*}
$\phi(z)$ has the following properties: \\
\indent (a) $\phi(z)$ is holomorphic in $\Omega$, \\
\indent (b) $|\phi(z)| \leq C$ in $\Omega$, \\
\indent (c) $|\phi(z)| \geq C^{-1} >0$ on $[-\epsilon/2, \epsilon/2]$, and \\
\indent (d) $|\phi(z)| \leq Ce^{-C/\delta}$ on $\Omega \cap \{ \Re z = \pm \epsilon \}$. \\
Now for $z \in \tilde{\Omega}:= [-\epsilon, \epsilon] + i[- \delta, \delta/ \log M]$ set  
\begin{eqnarray*}
g(z) = e^{-i N z \log M / \delta} \phi(z) f(z),
\end{eqnarray*}
and note that $g(z)$ satisfies \\
\indent (a) $|g(z)| \leq C M^{1-N}$ on $\tilde{\Omega} \cap \{ \Im z = - \delta \}$, \\
\indent (b) $|g(z)| \leq C_N e^{-C/\delta}$ on $\tilde{\Omega} \cap \{ \Re z = \pm \epsilon\}$, and \\
\indent (c) $|g(z)| \leq C_N \log(M) / \delta$ on $\tilde{\Omega} \cap \{ \Im z = \delta/ \log M \}$. \\
Then the classical maximum principle implies for $\delta$ sufficiently small and $N$ sufficiently large, 
$|g(z)| \leq C \log(M) / \delta$, which in turn implies
\begin{eqnarray*}
|f(z)| \leq C \frac{ \log M}{\delta} \,\,\, \text{on} \,\, \left[ -\frac{\epsilon}{2}, \frac{\epsilon}{2} \right] 
\subset \reals.
\end{eqnarray*}
\indent For part (b), note that our assumptions on $F(z)$ imply
\begin{eqnarray*}
\Im z \| u \|^2 \leq C \Re \langle F(z)u, u \rangle.
\end{eqnarray*}
We have 
\begin{eqnarray*} \|B F^{-1} A \|_{L^2 \to L^2} = \sup_{\{\|b\|_{L^2}
    = 1\}} \|B F^{-1} A b \|_{L^2} = \sup \| BF^{-1} A^2 b\|_{L^2},
\end{eqnarray*}
since $A^2 = A$.  Suppose $F(z) u(x) = A b(x,z)$. 
Then $u = F(z)^{-1} A b$ and $Bu = B F^{-1} AA b$, and for $\Im z >0$, 
\begin{eqnarray*}
\|Bu\|^2 & \leq & C \| u\|^2 \\
& \leq & \frac{C}{\Im z} \langle \Re F(z) u, u \rangle \\
& \leq & \frac{C}{\Im z} \left| \langle F(z) u, u \rangle \right|\\
& = &  \frac{C}{\Im z} \left| \langle Ab, u \rangle \right|\\
& = & \frac{C}{\Im z} \left| \langle Ab, Au \rangle \right| \\
& \leq & \frac{C}{\Im z} \| Ab\|^2
\end{eqnarray*}
where we have used $A^*A = A^2 = A$.  Thus we have 
\begin{eqnarray*}
 \|B F(z)^{-1}A\|_{L^2 \to L^2} & \leq & \frac{C}{ \sqrt{ \Im z}}, \,\,\, \text{for} \,\, \Im z >0 \,\, \text{and} \\
 \|B F(z)^{-1}A\|_{L^2 \to L^2} & = & \sup_{\{ \| u \| = 1 \}} \|
  BF^{-1}A u \|_{L^2} \\ 
& = & \sup_{\{ \| u \| = 1 \} } \| BF^{-1} BA u \|_{L^2
    \to L^2} \\
& \leq & M \sup_{\{ \| u \| = 1 \} } \| A u \|_{L^2} \\
& \leq & CM,
\end{eqnarray*}
and we can apply the proof of part (a) to $f(z) = B F(z)^{-1} A$ to get (b).
\end{proof}
\begin{proof}[Proof of Theorem \ref{main-theorem-2}]
Let $\psi_0$ satisfy the assumptions of Proposition \ref{Q(z)u-prop}.  Then
\begin{eqnarray*}
\| (1 -\psi_0)^w u \| \leq Ch^{-1} \| Q(z) u\| + \O(h^\infty) \|u\|.
\end{eqnarray*}
Further, since 
\begin{eqnarray*}
\left\| \left[ Q, \psi_0^w \right] u \right\| \leq \left\| \left[ Q, \psi_0^w \right] (1 - \tilde{\psi}_0^w)u \right\| 
+ \O(h^\infty) \|u \|,
\end{eqnarray*}
for some $\tilde{\psi}_0$ satisfying the assumptions of Proposition \ref{Q(z)u-prop} and $\WF \tilde{\psi}_0 \subset 
\{ \psi_0 = 1\}$, so using Theorem \ref{main-theorem-1} and the fact that $[Q, \psi_0^w]$ is compactly supported and 
of order $h$, we have
\begin{eqnarray*}
\| \psi_0^w u \| & \leq & Ch^{-N_0} \left( \| \psi_0^w Q u \| + \left\| \left[ Q, \psi_0^w \right] u \right\| \right) \\
& \leq & C h^{-N_0} \left(  \|\psi_0^w Q u\| + h^{-1} \| h Q u\| \right) + \O(h^\infty) \\
& \leq & C h^{-N_0} \| Qu \| + \O(h^\infty)\| u \|.
\end{eqnarray*}
\indent Now let $F(w)$ be the family of operators $F(w) = ih^{-1} Q(z_0 + hw)$, $A = \chi_{\supp \phi}^w$, $B = \id$. 
 Fix $\delta >0$ independent of $h$, $\epsilon = (Ch)^{-1}$, $M = h^{-N_0}$, and apply Lemma \ref{BZ-lemma} to get
\begin{eqnarray*}
\|B F^{-1} B \| & \leq & C \log (h^{-N_0}); \\
\|B F^{-1} A \| & \leq & C \sqrt{ \log (h^{-N_0}) },
\end{eqnarray*}
and (\ref{main-theorem-2-est-1}-\ref{main-theorem-2-est-2}) follows.
\end{proof}
\subsection{Proof of the Main Theorem}
The Main Theorem is an easy consequence of Theorem \ref{main-theorem-2}.
\begin{proof}[Proof of the Main Theorem]
Recall $A$ is $0$ microlocally away from $\gamma \times \gamma$.  Let $\tilde{A} \in \Psi_h^{0,0}$ be a pseudodifferential 
operator so that $\tilde{A} = I$ microlocally on a neighbourhood of
$\WF (A) \times \WF (A)$.  Let $a^w$ be as in 
Theorems \ref{main-theorem-1} and \ref{main-theorem-2}.  Choosing $A$ and $\tilde{A}$ so that $\WF (a^w)$ is 
disjoint from $\WF(\tilde{A})$, we have for $Q = Q(0)$
\begin{eqnarray}
\label{QAu}
Q\tilde{A} u = P \tilde{A} u.
\end{eqnarray}
The right hand side of \eqref{QAu} is $[ P, \tilde{A}] u + \tilde{A} Pu$.  Now $[P, \tilde{A}]$ is 
supported away from $\gamma$ since $\tilde{A}$ is constant near $\gamma$, so 
\begin{eqnarray}
\left\| P\tilde{A} u \right\|_{L^2(X)}  & \leq &  \left\| \left[ P,
  \tilde{A} \right] u \right\|_{L^2(X)}  + \|Pu\|_{L^2(X)} \nonumber \\ 
& \leq &  C h \left\| 
(I - A) u \right\|_{L^2(X)} + \| Pu \|_{L^2(X)}. \label{QAu-2}
\end{eqnarray}
From Theorem \ref{main-theorem-2}, we have
\begin{eqnarray}
\label{QAu-3}
\left\| Q \tilde{A} u \right\|_{L^2(X)} \geq C^{-1} \frac{h}{\sqrt{\log (1/h)} } \left\| \tilde{A} u \right\|_{L^2(X)}.
\end{eqnarray}
Combining \eqref{QAu-2} and \eqref{QAu-3}, we have
\begin{eqnarray*}
C^{-1} \|u \|_{L^2(X)} & \leq & C^{-1} \left( \left\|\tilde{A} u \right\|_{L^2(X)} + \left\| (I - A) u \right\|_{L^2(X)} 
\right) \\
& \leq & C \left( \sqrt{\log (1/h)} + C^{-1} \right) \left\| (I - A ) u
\right\|_{L^2(X)} \\
&& \quad \quad \quad \quad + C \frac{ \sqrt{ \log(1/h)}}{h} \| P u \|_{L^2(X)} ,
\end{eqnarray*}
which for $0 <  h < h_0$ is \eqref{main-theorem-3-est}.
\end{proof}
\begin{remark}
In the calculation \eqref{QAu-2}, we have only used $\|[P, \tilde{A}]u \| \leq Ch \| (I-A)u\|$.  If we could determine 
a global geometric condition which would allow us to choose $\tilde{A}$ in a manner which improves this, but doesn't have 
too much interaction with $a^w$ in the definition of $Q(z)$, we could eliminate the $\log(h^{-1})$ in 
\eqref{main-theorem-3-est}.
\end{remark}


\section{An Application: The Damped Wave Equation}
\label{wave-section}
\begin{figure}
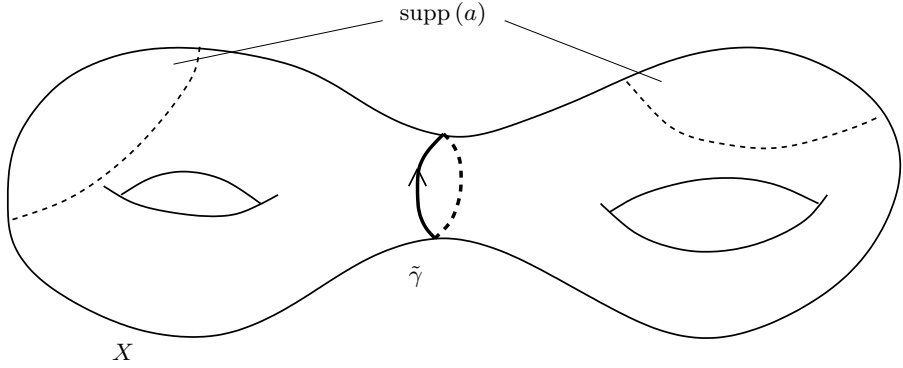

\include{fig5}
\caption{\label{figure:x-gamma} The manifold $X$ and the projection
  $\tilde{\gamma}$ of $\gamma$ onto $X$.}
\end{figure}

\numberwithin{equation}{section}
In this section we adapt the techniques from \S
\ref{main-theorem-1-proof}-\ref{main-theorem-2-proof} 
to study the damped wave equation.  Let $X$ be a compact manifold
without boundary, $a(x) \in \Ci(X)$, $a(x) \geq 0$, 
and consider the following problem:  
\begin{eqnarray}
\label{wave-equation-1}
\left\{ \begin{array}{l}
\left( \partial_t^2 - \Delta + 2a(x) \partial_t \right) u(x,t) = 0, \quad (x,t) \in X \times (0, \infty) \\
u(x,0) = 0, \quad \partial_t u(x,0) = f(x).
\end{array} \right.
\end{eqnarray}
Let $p \in \Ci (T^*X)$, 
$p= |\xi|^2$, be the microlocal principal symbol of $-\Delta$ and suppose the classical flow 
(geodesic flow) of $H_p$ admits a single closed, loxodromic orbit $\gamma$ in the 
level set $\{ p^{-1}(1)\}$.  Assume throughout that $a(x)$ is supported away from the 
projection $\tilde{\gamma}$ of $\gamma$ onto $X$ (see Figure \ref{figure:x-gamma}).  
We recall that the $H^s$ inner product on $X$ is given by the local formula
\be
\lll u, v \rrr_{H^s} = \int_{\reals^n} (1 + |\xi|^2)^s \hat{u}
\bar{\hat{v}} d\xi,
\ee 
where $\hat{u}$ is the Fourier transform of $u$.
If $u$ solves \eqref{wave-equation-1}, we define the {\it $s$-energy}
$E^s(t)$ of $u$ at time $t$ to be 
\begin{eqnarray*}
E^s(t) = \half \left( \left\| \partial_t u \right\|^2_{H^s(X)} + \left\|
\sqrt{-\Delta} u \right\|^2_{H^s(X)} \right).
\end{eqnarray*}

\begin{lemma}
If $a(x) \equiv 0$, $E^s(t)$ is constant.  If $a(x)$ is not identically zero, then $E^s(t)$ is 
decreasing.
\end{lemma}
\begin{proof}
\be
\frac{d}{dt} E^s(t) & = & \lll \partial_t^2 u, \partial_t u \rrr_{H^s} + \lll
\partial_t \sqrt{-\Delta}u, \sqrt{-\Delta} u \rrr_{H^s} \\
& = & \lll \partial_t u , ( \partial_t^2 - \Delta ) u \rrr \\
& = & -\lll \partial_t u, 2 a(x) \partial_t u \rrr.
\ee
\end{proof}

We make an important dynamical assumption, which amounts to a
geometric control condition similar to that given by
Rauch-Taylor in \cite{RT}.  We assume:
\begin{eqnarray}
\label{RT}
\left\{ \begin{array}{l} \text{There exists a time }T>0 \text{ and a
    neighbourhood } V \\
  \text{of } \gamma \text{ such that for all } |\xi| = 1,\,\,
  (x, \xi) \in T^*X \setminus V,\\
  \exp(tH_p)(x, \xi) \cap \{ a>0\}
    \neq \emptyset 
    \text{ for some } |t| \leq T.
    \end{array} \right.
\end{eqnarray}

In \cite{EvZw} \S $5.3$, it is shown that with a global Rauch-Taylor condition,
we have exponential decay in zero-energy.  Here we have a region
without geometric control, so we expect some loss.
\begin{theorem}
\label{exp-decay}
Assume \eqref{RT} holds and $a(x)$ is not identically zero.  Then for
any $\epsilon>0$, there is a constant $C>0$ such that 
\begin{eqnarray*}
E^0(t) \leq C e^{ -t/C } \| f\|_{H^\epsilon}^2.
\end{eqnarray*}
\end{theorem}

The damped wave equation in the context of a global Rauch-Taylor
condition has been studied in
\cite{RT}, \cite {Sjo2}, \cite{Leb}, and \cite{Hit}.  The difference
here is the presence of $\gamma$ and a neighbourhood in which the
Rauch-Taylor condition doesn't hold.

Formally, if $u \equiv 0$ for $t<0$, we apply the Fourier transform to \eqref{wave-equation-1} in
the $t$ variable and integrating by parts motivates us to study the
equation 
\begin{eqnarray}
\label{wave-equation-2}
P(\tau)\hat{u}(x, \tau) :=(-\tau^2 - \Delta + 2ia(x) \tau) \hat{u}(x,
\tau) = f.
\end{eqnarray}

We use the techniques of the previous sections to gain estimates on
the resolvent $P(\tau)^{-1}$.  We call the poles of $P(\tau)^{-1}$
{\it eigenfrequencies} for \eqref{wave-equation-1}.  Note if $\tau$ is
an eigenfrequency, then $0 \leq \Im \tau \leq 2 \|a\|_{L^\infty}$.
Further, \eqref{wave-equation-2} is invariant under the transformation
$(\hat{u}, \tau) \mapsto (\bar{\hat{u}}, - \bar{\tau})$, so the set of
eigenfrequencies is symmetric about the imaginary axis.  We
therefore study only those in the right half plane.  For $0 < h \leq
h_0$ and $z \in
\Omega: = [\alpha, \beta] + i[-\gamma, \gamma]$ where $0 <
\alpha<1<\beta<\infty$ and $\gamma>0$, set $\tau = \sqrt{z}/h$.
\eqref{wave-equation-2} becomes 
\begin{eqnarray}
\label{wave-equation-3}
\frac{1}{h^2} Q(z, h)\hat{u} = f
\end{eqnarray}
where 
\begin{eqnarray}
\label{Q(z,h)}
Q(z,h) = P(h) - z + 2ih \sqrt{z} a(x)
\end{eqnarray} 
and the principal symbol
of $P(h)$ is $p(x, \xi) = |\xi|^2$.  The next Corollary follows
directly from the proof of Theorem \ref{main-theorem-1}, replacing $s$
in the conjugation \eqref{Pth1} with $-s$.
\begin{corollary}
Suppose $u$ has wavefront set sufficiently close to $\gamma$.  Then
there exists $c_0 >0$, $C < \infty$, and $N \geq 0$ such that 
for $z \in [\alpha, \beta] + i[-c_0h, c_0h]$, 
\begin{eqnarray*}
\Im \left\langle (P(h) - z) u, u \right\rangle \geq C^{-1} h^N
\|u\|^2.
\end{eqnarray*}
In particular, $\|Q(z,h) u\| \geq C^{-1} h^N \|u\|$.
\end{corollary}

We observe that for $u$ as in the theorem and $-c_0 h <\Im z <0$, $\|Q(z, h) u \| \geq C^{-1} \Im z \|u\|$.  

The proof of Proposition
\ref{Q(z)u-prop} relies on the assumption that the symbol $a(x,
\xi)$ in \eqref{Q(z)} is elliptic away from $\gamma$.  The function
$a(x)$ in \eqref{Q(z,h)} is not assumed to be elliptic anywhere, so we
will use a technique from \cite{Leb} to replace
$a(x)$ with its average over trajectories of $\exp(tH_p)$.

For $T >0$, we define the $T$-trajectory average of a smooth function
$b$:
\begin{eqnarray*}
\langle b \rangle_T (x, \xi) = \frac{1}{T} \int_0^T b \circ
\exp(tH_p)(x, \xi) dt.
\end{eqnarray*}
Set $q(z) = 2 \sqrt{z}a(x)$, and for $z \in \widetilde{\Omega}:= [\alpha, \beta] + i[0, c_1h]$, where
$c_1 >0$ will be chosen later, and $(x, \xi) \in \{p^{-1}([\alpha -
\delta, \beta+ \delta])\}$ for $\delta>0$, let $g_{\Re z} \in \s (1)$
depending on $T$ solve 
\begin{eqnarray*}
q( \Re z) - H_pg_{\Re z} = \langle q( \Re z) \rangle_T.
\end{eqnarray*}
(See \cite{Sjo2} for details on the construction of $g_{\Re z}$.)  Now we form the elliptic
operator $A: = \Op_h^w (e^g) \in \Psi^{0,0}$, and observe
\begin{eqnarray*}
A^{-1} P A & = & P + A^{-1}[P,A] \\
& = & P - ih \Op_h^w(e^g)^{-1} \Op_h^w (\{p,e^g\}) \\
& = & P - ihB,
\end{eqnarray*}
with $\sigma_h(B) = e^{-g} \{p, e^g \} + \O(h) = H_pg + \O(h)$.  Thus
\begin{eqnarray*}
A^{-1} \left( P + ih q(z) \right) A & = & P + ih \Op_h^w \left( q( \Re z)
- H_pg \right) + \O(h^2) \\
& = & P + ih\Op_h^w \left( \langle q( \Re z ) \rangle_T \right),
\end{eqnarray*}
since $\Im z = \O(h) \Re z$.  Following \cite{Hit}, we claim there
exists a time $T>0$ such that 
\begin{eqnarray}
\label{inf-ave}
\langle a \rangle_T (x, \xi) \geq C^{-1} >0
\end{eqnarray}
for $(x, \xi) \in \{ p^{-1}([\alpha - \delta/2, \beta + \delta/2])
\}\setminus V$, where $V$ is as in the statement of Theorem \ref{exp-decay}.  To see this, recall $p = | \xi|^2$ means $H_p = 2 \langle \xi,
\partial_x \rangle$ and $p^{-1}(E) = \{ | \xi| = \sqrt{E} \}$, which
means 
\begin{eqnarray*}
\inf_{p^{-1}(E)} \langle a \rangle_T = \inf_{p^{-1}(1)} \langle a 
\rangle_{\sqrt{E}T}.
\end{eqnarray*}
By Assumption \eqref{RT}, 
\begin{eqnarray*}
\inf_{p^{-1}(1)} \langle a 
\rangle_{\sqrt{E}T} \geq C^{-1} >0
\end{eqnarray*}
in $\{p^{-1}(1)\} \setminus V$ for $T$ sufficiently large and
$\sqrt{E}$ close to $1$.  Choosing
$\alpha$ and $\beta$ sufficiently close to $1$ proves \eqref{inf-ave}.
\begin{corollary}
\label{Q(z)u-prop-2}
Suppose $\psi_0 \in \s^{0,0}(T^*X) \cap \Ci_c(T^*X)$ is a microlocal cutoff function to a small neighbourhood 
of $\gamma \subset \{p^{-1}(1) \}$.  For $Q(z,h) = P(h) - z
+2ih\sqrt{z}a$ as above with $z \in [\alpha, \beta] + i(-c_1h, c_1h)$, 
$c_1 >0$, we have
\begin{eqnarray}
\label{Q(z)u-prop-est-2}
Q(z,h)u = f \Longrightarrow \left\| (1 - \psi_0 )^w u \right\| \leq C h^{-1} \|f \| + \O(h^\infty) \|u \|.
\end{eqnarray}
\end{corollary}
\begin{proof}
Selecting $T >0$ sufficiently large and $c_1 >0$ such that 
\begin{eqnarray*}
0 < c_1 < \inf_{p^{-1}([\alpha - \delta/2, \beta + \delta/2])} \langle a \rangle_T (x, \xi),
\end{eqnarray*}
we apply the proof of Proposition \ref{Q(z)u-prop} to the conjugated
operator $A^{-1} Q(z,h)A$.
\end{proof}
We now have good resolvent estimates for $z$ in an $h$ interval below the
real axis, as well as weaker estimates above.
\begin{corollary}
(i) There exist constants $C>0$ and $N >0$ such that the resolvent $Q(z,h)^{-1}$ satisfies
\begin{eqnarray*}
 \|Q(z,h)^{-1}\|_{L^2 \to L^2} \leq C h^{-N} , \quad z \in
   [\alpha, \beta] + i(-c_0h, c_0h).
\end{eqnarray*}
(ii) In addition, there is a 
constant $C_1$ such that
\begin{eqnarray*}
\|Q(z,h)^{-1}\|_{L^2 \to L^2} \leq C_1 \frac{\log(1/h)}{h}, \quad z
\in [\alpha, \beta] +i[-C_1^{-1}h / \log(1/h), C_1^{-1}h].
\end{eqnarray*}
\end{corollary}
This is an immediate consequence of the proof of Theorem
\ref{main-theorem-2}, together with the slight modification of Lemma
\ref{BZ-lemma} given in Lemma \ref{holo-lemma-2}.  

\begin{lemma}
\label{holo-lemma-2}
Let $f(z)$ be a holomorphic function on $\Omega = [-\epsilon,
  \epsilon] + i [-\delta, \delta]$, with
\be
\frac{\delta}{\epsilon} \ll M^{-\frac{1}{N_1}}
\ee
for some $N_1>0$, and suppose 
$f$ satisfies $|f(z) | \leq M$ on
$\Omega$ with $|f(z)| \leq  C|\Im z|$ for $\Im z <0$.  Then there
exists a constant $0<C_1<\infty$ such that if $ - C_1^{-1}\delta/\log M \leq \Im z
\leq C_1^{-1}\delta$ we have
\be
|f(z)| \leq C \frac{\log M}{\delta}.
\ee
\end{lemma}
\begin{proof}
Let $\psi(x)$ be as in the proof of Lemma \ref{BZ-lemma}, and for
$C_1^{-1}\ll c_0$, let
\be
\phi(z) = \delta^{-\half} \int e^{-(x -z + iC_1^{-1} \delta )^2/ \delta} \psi(x) dx.
\ee
We observe if $C_1>0$ is sufficiently large, for $|\Im z| \leq
C_1^{-1} \delta $, 
\be
\lefteqn{(x - z + i C_1^{-1} \delta)^2  = } \\ & = & (x - \Re z)^2 - (C_1^{-1}\delta -
\Im z)^2 + 2i(x - \Re z)(C_1^{-1}\delta - \Im z) 
\ee
and
\be
\left| (x - \Re z) ( C_1^{-1} \delta - \Im z) \right| \leq 4 C_1^{-1}
\epsilon \delta,
\ee
so if $C_1>0$ is sufficiently large,
\be
\Re e^{-(x -z + iC_1^{-1} \delta )^2/ \delta} & \geq & e^{-(x-\Re
  z)^2/\delta + (C_1^{-1} \delta - \Im z )^2/\delta} \cos(4C_1^{-1}
\epsilon) \\
& \geq & C^{-1}
e^{-(x-\Re z)^2/\delta + (C_1^{-1} \delta - \Im z )^2/\delta}.
\ee
Thus$\phi(z)$ satisfies

(a) $\phi(z)$ is holomorphic in $\Omega$,

(b) $|\phi(z)| \leq C$ in $\Omega$,

(c) $|\phi(z)| \geq C^{-1}$ for $z \in [-\epsilon/2, \epsilon/2] +
i[-C_1^{-1}\delta, C_1^{-1} \delta]$,

(d) $|\phi(z)| \leq C e^{-C/\delta}$ on $\{ \pm \epsilon \} \times
i[-C_1^{-1} \delta, C_1^{-1} \delta]$, if $C_1>0$ is chosen large enough.

Now similar to the proof of Lemma \ref{BZ-lemma}, for 
\be
z \in
\tilde{\Omega}:= [-\epsilon, \epsilon] + i[- C_1^{-1} \delta/ \log M,
  C_1^{-1}\delta]
\ee
set  
\begin{eqnarray*}
g(z) = e^{i N z \log M / \delta} \phi(z) f(z).
\end{eqnarray*}

Then as in the proof of Lemma \ref{BZ-lemma}, the classical maximum principle implies for $\delta$ sufficiently small and $N$ sufficiently large, 
$|g(z)| \leq C \log(M) / \delta$, which in turn implies
\begin{eqnarray*}
|f(z)| \leq C \frac{ \log M}{\delta} \,\,\, \text{on} \,\, \left[
  -\frac{\epsilon}{2}, \frac{\epsilon}{2} \right] + i [-C_1^{-1}\delta / \log
  M , C_1^{-1} \delta ].
\end{eqnarray*}
\end{proof}

With these resolvent estimates, we have the following estimates in
terms of $\tau$.
\begin{proposition}
\label{P-tau-bnds}
Fix $\epsilon>0$.  There exist constants $0 < C, C_1< \infty $ such that if
\be
-(\log \lll \tau \rrr)^{-1} \leq \Im \tau \leq C_1^{-1},
\ee
\begin{eqnarray}
\label{P-L2} \| P( \tau)^{-1} \|_{L^2 \to L^2} & \leq &
\frac{C\log \lll \tau \rrr}{\langle \tau \rangle }, \\
\label{P-H2}  \| P( \tau)^{-1} \|_{L^2 \to H^2} & \leq &
{C}{\langle \tau \rangle \log \lll \tau \rrr}, \,\, \text{and} \\
\label{P-H1}   \| P( \tau)^{-1} \|_{H^s \to H^{s+1-\epsilon}} & \leq & C.
\end{eqnarray}
\end{proposition}
\begin{proof}
\eqref{P-L2} follows directly from rescaling.  To see \eqref{P-H2},
calculate
\begin{eqnarray*}
\|u \|_{H^2} & \leq & C  ( \| \Delta u \|_{L^2} + \| u \|_{L^2} ) \\
& \leq & C \left( \left\| P(\tau) u \right\|_{L^2} + \left\| ( -
\tau^2 + 2ia(x) \tau ) u \right\|_{L^2} + \frac{\log \langle \tau \rangle}{\langle \tau
  \rangle } \left\| P(
\tau) u \right\|_{L^2} \right) \\
& \leq & C \left(1 + | \tau |\log \lll \tau \rrr + \langle \tau
\rangle^{-1} \log \lll \tau \rrr \ \right)
\left\| P( \tau) u \right\|_{L^2}.
\end{eqnarray*}
For \eqref{P-H1}, let $\epsilon>0$ be given.  From Lemma
\ref{inter-lemma}, we have
\be
\left\| P(\tau)^{-1} u \right\|_{H^{1-\epsilon}}^2 & \leq & C \left \|
P(\tau)^{-1} u
\right\|_{H^2}^{1-\epsilon} \left\| P(\tau)^{-1}u \right\|_{L^2}^{1 +
  \epsilon} \\
& \leq & C_\epsilon \|u\|_{L^2(X)}^2.
\ee
To get the estimates for $H^s \to H^{s+1-\epsilon}$, we conjugate
$P(\tau)^{-1}$ by the operators 
\be
\Lambda^s = (1 - \Delta)^\frac{s}{2}
\ee
and apply to $v = \Lambda^s u$:
\be
\left\|P(\tau)^{-1} u \right\|_{H^{s+1-\epsilon}} & = &
\left\|\Lambda^{1-\epsilon} \Lambda^sP(\tau)^{-1} \Lambda^{-s} v
\right\|_{L^2} \\
& = & \left\| \Lambda^{1-\epsilon} \left( P(\tau)^{-1}
+ \Lambda^s [ P(\tau)^{-1} , \Lambda^{-s}] \right) v \right\|_{L^2} \\
& \leq & \left\| v \right\|_{L^2}  \\
& \leq & C \left\|u \right\|_{H^s}.
\ee
\end{proof} 
We have used the following interpolation lemma.
\begin{lemma}
\label{inter-lemma}
Let $\epsilon>0$ be given, and suppose $f \in H^2(X) \cap
L^2(X)$.  Then
\be
\left\|f\right\|_{H^{1-\epsilon}}^2 \leq C \left \| f
\right\|_{H^2}^{1-\epsilon} \left\| f \right\|_{L^2}^{1 + \epsilon}.
\ee
\end{lemma}
\begin{proof}
We use the local formula for $H^s$ norms and calculate:
\be
\left\|f\right\|_{H^{1-\epsilon}}^2 &=& \int_{\reals^n} (1 +
|\xi|^2)^{1-\epsilon} \hat{f} \bar{\hat{f}} d \xi \\
& = & \int \left( (1 + |\xi|^2) |\hat{f}| \right)^{1-\epsilon}
|\hat{f}|^{1+\epsilon} d \xi \\
& \leq & C \left\|\left( (1 + |\xi|^2) |\hat{f}| \right)^{1-\epsilon}
\right\|_{L^{\frac{2}{1-\epsilon}}} \left\| |\hat{f}|^{1+\epsilon}
\right\|_{L^{\frac{2}{1+\epsilon}}} \\
& = & \leq C \left \| f
\right\|_{H^2}^{1-\epsilon} \left\| f \right\|_{L^2}^{1 + \epsilon} .
\ee
\end{proof}

We are now in position to prove Theorem \ref{exp-decay}.  This proof
comes almost directly from \cite{EvZw} \S $5.3$.
\begin{proof}[Proof of Theorem \ref{exp-decay}]
Assume $u(x,t)$ solves \eqref{wave-equation-1}.  
Choose $\chi \in \Ci( \reals )$, $0 \leq \chi \leq 1$, $\chi \equiv 1$
on $[1, \infty)$, and $\chi \equiv 0$ on $(-\infty, 0]$.  Set $u_1 (x,t)=
\chi(t) u(x,t)$.  We apply the damped wave operator to $u_1$:
\begin{eqnarray}
\lefteqn{ \left( \partial_t^2 - \Delta + 2a \partial_t \right) u_1 =}
\label{wave-equation4.a} \\ & = & \chi''
u + 2 \chi' u_t + 2a \chi' u + \chi \left( \partial_t^2 - \Delta + 2a \partial_t
\right) u \nonumber \\
& = & \chi''
u + 2 \chi' u_t + 2a \chi' u =: g_1. \label{wave-equation4.b}
\end{eqnarray}
With $g_1$ supported in $X \times (0,1)$ and $u_1 \equiv 0$ for $t\leq
0$, we have
\begin{eqnarray}
\label{g1-bnd}
\|g_1\|_{L^2((0, \infty); H^\epsilon)}^2 \leq C \left( \| u \|_{L^2((0,1);
  H^\epsilon)}^2 + \|\partial_t u\|_{L^2((0,1);H^\epsilon)}^2 \right).
\end{eqnarray}
Now
\begin{eqnarray*}
\partial_t \langle u,u \rangle_{H^\epsilon(X)} & = & 2 \langle \partial_t u, u
\rangle_{H^\epsilon(X)}  \\
& \leq & \| \partial_t u \|_{H^\epsilon(X)}^2 + \| u \|_{H^\epsilon(X)}^2 \\
& \leq & CE^\epsilon(t) + \|u \|_{H^\epsilon(X)}^2,
\end{eqnarray*}
so by Gronwall's inequality, 
\begin{eqnarray*}
\|u(t, \cdot) \|_{H^\epsilon(X)}^2 & \leq& C e^t\left( \|u(0,\cdot)
\|_{H^\epsilon(X)}^2 + \int_0^t E^\epsilon(s)ds \right) \\
& \leq & Ct e^t \| f \|_{H^\epsilon(X)}^2.
\end{eqnarray*}
Thus \eqref{g1-bnd} is bounded by $C \|f \|_{H^\epsilon(X)}^2$.

We now apply the Fourier transform to
(\ref{wave-equation4.a}-\ref{wave-equation4.b}) to write $\hat{u}_1 =
P(\tau)^{-1} \hat{g}_1$.  By Proposition \ref{P-tau-bnds}, we have for $\Im \tau =
C^{-1}>0 $
\begin{eqnarray*}
\left\| e^{t/C} u_1 \right\|_{L^2((0, \infty); H^1)} & = & \left\|
\hat{u}_1( \cdot + iC^{-1}) \right\|_{L^2((-\infty, \infty);H^1)} \\
& = & \left\| P(\cdot + iC^{-1})^{-1} \hat{g}_1( \cdot + i C^{-1})
\right\|_{L^2((-\infty, \infty); H^1)} \\
& \leq & C \|\hat{g}_1 \|_{L^2((-\infty, \infty); H^\epsilon)} \\ 
& \leq & C\| g_1\|_{L^2((0, \infty); H^\epsilon)} \\
& \leq & C\|f\|_{H^\epsilon(X)}.
\end{eqnarray*}
Thus
\begin{eqnarray*}
\|e^{t/C} u \|_{L^2((1, \infty);H^1)} \leq C \|f\|_{H^\epsilon(X)}.
\end{eqnarray*}
Now for $T>2$, choose $\chi_2 \in \Ci( \reals)$, $0 \leq \chi_2 \leq 1$,
such that $\chi_2 \equiv 0$ for $t \leq T-1$, and $\chi_2 \equiv 1$ for $t
\geq T$.  Set $u_2 (x,t) = \chi_2(t) u(x,t)$.  We have
\begin{eqnarray*}
\left( \partial_t^2 - \Delta + 2a \partial_t \right) u_2 = g_2
\end{eqnarray*}
for $g_2 = \chi_2'' u + 2 \chi_2' u_t + 2a \chi_2' u$, and $\supp g_2
\subset X \times [T-1,T]$.  Define
\begin{eqnarray*}
E_2(t) = \frac{1}{2} \int_X \left( \partial_t u_2 \right)^2 + \left|
\sqrt{-\Delta} u_2 \right|^2 dx,
\end{eqnarray*}
and observe
\begin{eqnarray*}
E_2'(t) & = & \left\langle \partial_t^2 u_2, \partial_t u_2
\right\rangle_X - \left\langle \Delta u_2, \partial_t u_2
\right\rangle_X \\
& = & - \left\langle 2a(x) \partial_t u_2, \partial_t u_2
\right\rangle_X + \left\langle g_2, \partial_t u_2 \right\rangle_X \\
& \leq & C \int_X \left| \partial_t u_2 \right| \left( \left|
\partial_t u \right| + |u| \right) \\
& \leq & C \left( E_2(t) + \int_X \left( \left| \partial_t u \right|^2 + |u|^2
\right) dx\right).
\end{eqnarray*}
Now since $E_2(T-1) = 0$ and $E_2(T) = E(T)$, Gronwall's inequality
gives
\begin{eqnarray}
\label{ET-1}
E(T) \leq C \left( \left\| \partial_t u \right\|_{L^2((T-1,T);L^2)}^2 +
\| u \|_{L^2((T-1,T); L^2)}^2 \right).
\end{eqnarray}
We need to bound the first term on the right hand side of
\eqref{ET-1}.  Choose $\chi_3 \in \Ci( \reals)$ such that $\chi_3
\equiv0$ for $t \leq T-2$ and $t \geq T+1$, $\chi_3 \equiv 1$ for $T-1
\leq t \leq T$.  Then
\begin{eqnarray*}
0 & = & \int_{T-2}^{T+1} \int_X \chi_3^2 u \left( \partial_t^2 u -
\Delta u + 2a \partial_t u \right) dxdt \\
& = & \int_{T-2}^{T+1} \int_X -\chi_3^2 (\partial_t u)^2 - 2 \chi_3
\chi_3' u \partial_t u + 2 \chi_3^2a \partial_t u + \chi_3^2 |
\sqrt{-\Delta} u |^2 dxdt,
\end{eqnarray*}
whence
\begin{eqnarray*}
\left\| \partial_t u \right\|_{L^2((T-1, T);L^2)} \leq C \| u
\|_{L^2((T-2, T+1); H^1)},
\end{eqnarray*}
giving
\begin{eqnarray*}
E(T) \leq C \|u \|_{L^2((T-2, T+1);H^1)}^2 \leq C e^{-T/C}
\|f\|_{H^\epsilon(X)}^2
\end{eqnarray*}
as claimed.
\end{proof}


\end{document}